\documentclass[11pt]{amsart}
\usepackage{amsmath,amssymb,latexsym,amsfonts,amscd}
\newtheorem{theorem}{Theorem}[section]

\newtheorem{lemma}[theorem]{Lemma}
\newtheorem{corollary}[theorem]{Corollary}
\newtheorem{example}[theorem]{Example}

\newtheorem{remark}[theorem]{Remark}

\begin{document}

\title[on some local cohomology modules]{On some local cohomology
modules}
\author{Gennady Lyubeznik}
\address{Department of Mathematics, University of Minnesota, Minneapolis,
MN 55455}
\email{gennady@math.umn.edu}

\maketitle

\section{Introduction}
All rings in this paper are commutative and Noetherian. 
If $R$ is a ring and $I\subset R$ is an ideal, cd$(R,I)$ denotes
the cohomological dimension of $I$ in $R$, i.e. the largest integer $i$
such that the $i$-th local cohomology
module $H^i_I(M)$ doesn't vanish for some $R$-module $M$. 

For the purposes of this introduction 
$R$ is a complete equicharacteristic regular local
$d$-dimensional ring with a separably closed residue field and $I\subset
R$ is an ideal such that every minimal prime over $I$ has height
at most $c$. We quote the following two results.

\medskip 

(i) cd$(R,I)\leq
d-[\frac{d-1}{c}]$ \cite[Korollar 2]{Fa} (see Theorem \ref{Falt} below).
This bound is sharp for all $d$ and $c$ \cite{Lyub}.

(ii) If,
in addition,
$I$ is prime and $c<d$, then cd$(R,I)\leq d-[\frac{d-2}{c}]-1$
\hbox{\cite[3.8]{HL}} (see Theorem \ref{HL} below). This is also
sharp for all $d$ and $c$ \cite[5.6]{HL}.

\medskip

These two upper
bounds, $v=d-[\frac{d-2}{c}]-1$ and $v'=d-[\frac{d-1}{c}]$, coincide
if
$c |(d-1)$, otherwise $v'=v+1$. This
raises the following natural question.

\medskip

\noindent{\bf Question.} {\it Under what
conditions is  ${\rm cd}(R,I)\leq
v$? Equivalently, under what 
conditions is $H_I^{v+1}(M)=0$
for every $M$?}

\medskip

We discussed this question in our old paper \cite{HL} in the special case
that $[\frac{d-2}{c}]=2$ \cite[5.8-5.12]{HL}, but were unable to settle
even this special case. The purpose
of the present paper is to provide a complete answer to this question.
Namely, we give a necessary and sufficient condition 
exclusively in terms of combinatorial properties of the set of subsets
$\{i_0,\dots,i_j\}$ of the indexing set $\{1,\dots,n\}$ of the minimal
prime ideals $I_1,\dots, I_n$ over $I$, such that 
the sum $I_{i_0}+\dots +I_{i_j}$ is
$\mathfrak m$-primary. More precisely, we prove the following.

\begin{theorem}\label{maintheorem}
Let $R$ be a complete $d$-dimensional equicharacteristic regular local
 ring with maximal ideal $\mathfrak m$ and a separably closed residue
field $k$. Let
$I\subset R$ be an ideal of $R$ and let $M$ be an
$R$-module. Let
$I_1,\dots,I_n$ be the minimal primes  of $I$. Assume the height of every
$I_i$ is at most $c<d$. Let
$\Delta$ be the simplicial complex on the vertices $1,\dots,n$ such
that a simplex
$\{i_0,\dots,i_s\}$ is included in $\Delta$ if and only if
$I_{i_0}+\dots+I_{i_s}$ is 
$\underline{not}$ $\mathfrak m$-primary. Let $t=[\frac{d-2}{c}]$ and let
$v=d-[\frac{d-2}{c}]-1$. Then 
$H^{v+1}_I(M)$ is isomorphic
to the direct sum of $w$ copies of $H^d_{\mathfrak m}(M)$, where 
$w={\rm dim}_k\tilde H_{t-1}(\Delta;k)$ and
 $\tilde H_*(\Delta;k)$ is the reduced singular homology of $\Delta$
with coefficients in $k$. In particular, 
\hbox{$cd(R,I) \leq v$} if and 
only if $\tilde H_{t-1}(\Delta;k)=0$.
\end{theorem}

If $c=d-2$, then $t=[\frac{d-2}{c}]=1$ and the theorem implies the
following: {\it if an ideal I is the intersection of prime ideals of
heights at most $d-2$, then ${\rm cd}(R,I)\leq d-2$ if and only if the
punctured spectrum of $R/I$ is connected} (indeed, the punctured
spectrum is connected iff $\Delta$ is connected, i.e. iff \hbox{$\tilde
H_0(\Delta;k)=0$).} This is a famous result with a rich history.
Various pieces of this result have been proven by Peskine and Szpiro
\cite[III, 5.5]{PSz}, Hartshorne \cite[7.5]{Ha}, Ogus \cite[2.11]{O} and 
Huneke and Lyubeznik \cite[2.9]{HL}. Theorem \ref{maintheorem} extends
this famous result to higher values of $t=[\frac{d-2}{c}]$.

Theorem \ref{maintheorem} makes it plain why the sharp bounds $v$ and $v'$
of (i) and (ii) coincide in the case that $c|(d-1)$. Indeed, in this case
$d=c(t+1)+1$, where $t=[\frac{d-2}{c}]$, hence the sum of at most $t+1$
primes of height at most $c$ each cannot add up to an $\mathfrak
m$-primary ideal, i.e. $\Delta$ contains every single simplex
$\{i_0,\dots,i_p\}$ of dimension $p\leq t$. Hence the reduced singular
homology of $\Delta$ in degrees at most $t-1$ coincides with that of the
full simplex on $\{1,\dots,n\}$. But the reduced singular homology of the
latter vanishes in all degrees, i.e. $\tilde H_{t-1}(\Delta;k)=0$ for
all $I$. Hence cd$(R,I)\leq v$ for all, not just prime, $I$.

Theorem \ref{maintheorem} immediately implies an extension of 
(ii) to the case that $I$ is not necessarily
prime, but has a small number of minimal primes. Namely, a sum
of  fewer than
$\frac{d}{c}$ primes of height at most $c$ cannot add up to an $\mathfrak
m$-primary ideal in a regular $d$-dimensional ring, hence if
$n<\frac{d}{c}$, then the complex $\Delta$ is the full simplex on
vertices
$\{1,\dots,n\}$ whose reduced singular homology vanishes in all degrees.
Thus we get the following corollary which is a generalization of (ii).

\begin{corollary}\label{maincoro}
Let $R$ be a complete $d$-dimensional equicharacteristic regular local
 ring with a separably closed residue field. Let
$I\subset R$ be an ideal with $n<\frac{d}{c}$ minimal
primes each of which has height at most $c$. Then \hbox{$cd(R,I)
\leq d - 1 - [\frac{d-2}{c}]$}.
\end{corollary}

The bound $n<\frac{d}{c}$ is sharp for all $d$ and $c$ such that
$c\not|(d-1)$ (if $c|(d-1)$, then cd$(R,I)\leq d - 1 - [\frac{d-2}{c}]$
for all $I$; see (i) above). Indeed, if $n$ is the smallest integer
$\geq
\frac{d}{c}$, then there exist prime ideals $I_1,\dots,I_n$ of height $c$
each, such that $I_1+\dots+I_n$ is $\mathfrak m$-primary in which case
$\Delta$ contains every simplex except the top-dimensional one, i.e.
$\Delta$ is homeomorphic to the $(n-2)$-sphere, hence $\tilde
H_{n-2}(\Delta;k)\ne 0$, i.e. cd$(R,I)>d-1-[\frac{d-2}{c}]$ by Theorem
\ref{maintheorem}, where $I=I_1\cap\dots\cap I_n$.

\medskip

The principal tool of our proofs is the following Mayer-Vietoris spectral
sequence of \cite[p. 39]{AGZ} (see Theorem \ref{MV} below).
$$E_1^{-p,q}=\oplus_{i_0<\dots<i_p}H^q_{I_{i_0}+\dots+I_{i_p}}(M)
\Longrightarrow H^{q-p}_{I_1\cap\dots\cap I_n}(M).$$ 
We use it in Section 2 to reduce the proofs of our main results to an
auxiliary statement that will be proven in Theorem \ref{main}, the main
theorem of Section 3 (see Corollaries \ref{plan} and \ref{MVDelta} and the
paragraph following the proof of 
\ref{MVDelta}). 

In the special case that
$[\frac{d-2}{c}]= 2$ that auxiliary statement readily follows from some
well-known results which gives a short proof of Theorem
\ref{maintheorem} in this special case. We conclude Section 2 by using
this special case of Theorem
\ref{maintheorem} to gain a complete understanding of an old
example that defied our efforts in \cite[5.12]{HL} (Example
\ref{exa}). The complex
$\Delta$ in this example turns out to be homeomorphic to the real
projective plane and consequently cd$(R,I)$ depends on whether
char$k=2$.

Section 3 is devoted to a proof of the above-mentioned Theorem
\ref{main}. 

The final section, Section 4,
contains proofs of our main results. We prove our main results for
arbitrary local rings containing a field. Theorem \ref{maintheorem} and
Corollary \ref{maincoro} are special cases of Corollaries \ref{Delta} and
\ref{Coro-gener} respectively.

\' Etale analogues of the results in this paper and their topological
applications will be dealt with in an upcoming paper \cite{L-upc}.

\section{The Mayer-Vietoris spectral sequence}\label{S:MV}
The Mayer-Vietoris spectral sequence of \cite[p. 39]{AGZ} is stated and
proven in \cite{AGZ} only for defining ideals of arrangements of
subspaces and only for
$M=A$. Yet the same statement\footnote
{There is a misprint in a key formula of \cite[p. 39]{AGZ} which says
$E_1^{-i,j}=Roos_i(H^0_{[*]}(E^j))$ whereas it should say
$E_1^{-i,j}=Roos_i(H^j_{[*]}(R))$.} holds and the same proof works in
general. We reproduce a complete proof to avoid any misunderstanding.
\begin{theorem}\label{MV}
Let $A$ be a commutative ring, let $I_1,\dots,I_n\subset A$ be ideals and
let $M$ be an $A$-module. There exists a spectral sequence
$$E_1^{-p,q}=\oplus_{i_0<\dots<i_p}H^q_{I_{i_0}+\dots+I_{i_p}}(M)
\Longrightarrow H^{q-p}_{I_1\cap\dots\cap I_n}(M).$$
\end{theorem}

\emph{Proof.} We recall that for an ideal $J$ of $A$ one denotes by
$\Gamma_J(M)$ the submodule of $M$ consisting of the elements
annihilated by some power of $J$. If $J'\subset J$, we let
$\gamma_{J,J'}:\Gamma_{J}(M)\hookrightarrow \Gamma_{J'}(M)$ be the natural
inclusion.

We denote by $\Gamma^{\bullet}(M)$ the complex 
$$0\to \Gamma^{-n+1}(M)\stackrel{d^{-n+1}}{\longrightarrow}
\Gamma^{-n+2}(M)\stackrel{d^{-n+2}}{\longrightarrow}
\dots\stackrel{d^{-1}}{\longrightarrow}
\Gamma^0(M)\to 0$$ in the category of $R$-modules where
$\Gamma^{-p}(M)=\oplus_{
i_0<i_1<\dots<i_p}\Gamma_{I_{i_0}+\dots+I_{i_p}}(M)$ and 
$d^{-p}(x)=\oplus_{j=0}^{j=p}(-1)^j\gamma_{J,
J_j}(x)$ for every element $x\in \Gamma_{J}(M)\subset\Gamma^{-p}(M)$,
where
$J=I_{i_0}+\dots+I_{i_p}$ and
$J_j=I_{i_0}+\dots+I_{j-1}+I_{j+1}+\dots+I_{i_p}$ and $\gamma_{J,
J_j}(x)$ is an element of $\Gamma_{J_j}(M)\subset \Gamma^{-p+1}(M)$. 

We view
$\Gamma(-)$ as a functor from the category of
$A$-modules to the category of complexes of $A$-modules; namely, an
$A$-module map $f:M\to M'$ iduces maps
$f_{i_0,\dots,i_p}:\Gamma_{I_{i_0}+\dots+I_{i_p}}(M)\to
\Gamma_{I_{i_0}+\dots+I_{i_p}}(M')$ which induce a map of complexes
$f^{\bullet}:\Gamma(M)\to \Gamma(M')$ where
$f^{-p}=\oplus_{
i_0<i_1<\dots<i_p}f_{i_0,\dots,i_p}$.

We claim that if $M$ is injective, then $H^t(\Gamma^{\bullet}(M))$, the
$t$-th cohomology module of the complex $\Gamma^{\bullet}(M)$, is zero
if $t<0$ and is isomorphic to $\Gamma_{I_1\cap\dots\cap I_n}(M)$ if $t=0$.
Indeed,  an injective module is a direct sum of modules of the form
$E(A/P)$ where $P$ is some prime ideal of $A$ and $E(A/P)$ is the
injective hull of $A/P$. Since the functors $\Gamma^{\bullet}(-)$ and 
$\Gamma_{I_1\cap\dots\cap I_n}(- )$ commute with direct sums, it is enough
to prove the claim for $M=E(A/P)$ 
in which case $\Gamma_J(M)=0$ if $J\not\subset P$ and
$\Gamma_J(M)=M$ if $J\subset P$ while the inclusion maps $\gamma_{J',J}$
for ideals $J\subset J'\subset P$ are the identity maps. If none of
the ideals $I_1,\dots,I_n$ are contained in $P$, then the complex
$\Gamma^{\bullet}(M)$ vanishes and so does the module
$\Gamma_{I_1\cap\dots\cap I_n}(M)$, hence the claim holds. If some
of the ideals $I_j$ are contained in $P$, then
$\Gamma^{\bullet}(M)=M\otimes_k\mathcal K^{\bullet}$ where $k$ is the
fraction field of $A/P$ and 
$\mathcal K^{\bullet}$ is the complex such that upon giving
$\mathcal K^{-p}$ the name
$\mathcal K_p$, the complex
$\mathcal K^{\bullet}$ turns into a complex $\mathcal K_{\bullet}$ which
is the standard algebraic complex for the computation of the singular
homology with coefficients in
$k$ of the full simplex on the set of vertices $\{i|I_i\subset
P\}$. Since a full simplex on any finite set of vertices in
contractible, the homology of $\mathcal K_{\bullet}$ vanishes in positive
degrees and its zeroth homology is
$k$. Hence the cohomology of $\mathcal K^{\bullet}$ vanishes in negative
degrees and its zeroth cohomology is $k$. Since tensoring over a
field is an exact functor and $M\otimes_kk=M=\Gamma_{I_1\cap\dots\cap
I_n}(M)$ for $M=E(A/P)$, the claim is proven.

 Since
$\Gamma(-)$ is a functor, any complex $C^{\bullet}$ yields a
double complex $\Gamma^{\bullet}(C^{\bullet})$ in which the $(-p,q)$-th
entry is $\Gamma^{-p}(C^q)$, the vertical map
$\Gamma^{-p}(C^q)\to \Gamma^{-p}(C^{q+1})$ is induced by the
differential $C^q\to C^{q+1}$ of $C^{\bullet}$ and the
$q$-th horizontal line is 
$\Gamma^{\bullet}(C^q)$ (with all the differentials multiplied by
$(-1)^q$, to make $\Gamma^{\bullet}(C^{\bullet})$ an honest double
complex, rather than a commutative diagram). Since the double complex 
$\Gamma^{\bullet}(C^{\bullet})$ has just a finite number of non-zero
columns, both associated spectral sequences converge to the cohomology of
the total complex of $\Gamma^{\bullet}(C^{\bullet})$.

Let $0\to M\to E^0\to E^1\to \dots$ be an injective resolution of $M$.
Setting $C^{\bullet}=\Gamma_{I_1\cap\dots\cap I_n}(E^{\bullet})$ 
we get a double complex
$\Gamma^{\bullet}(C^{\bullet})$ as above in which the \hbox{$q$-th}
horizontal line is
$\Gamma^{\bullet}(\Gamma_{I_1\cap\dots\cap I_n}(E^q))$ with all the
differentials multiplied by $(-1)^q$. Since
$E^q$ is injective, so is $\Gamma_{I_1\cap\dots\cap I_n}(E^q)$, hence
according to the above claim, the
$q$-th horizontal line is exact except in degree zero where the cohomology
is
$\Gamma_{I_1\cap\dots\cap I_n}(E^q)$. Hence the horizontal cohomology of
the double complex is concentrated on the vertical line $p=0$ where it
forms the complex $\Gamma_{I_1\cap\dots\cap I_n}(E^{\bullet})$ whose
cohomology in degree $q$ is $H^q_{I_1\cap\dots\cap I_n}(M)$. Thus
the abutment in total degree $t$ of each of the two spectral sequences
associated with $\Gamma^{\bullet}(C^{\bullet})$ is
$H^t_{I_1\cap\dots\cap I_n}(M)$. 

The
$E_1^{-p,q}$ term of one of the two associated spectral sequences
is the vertical cohomology of
$\Gamma^{\bullet}(C^{\bullet})$ in degree $(-p,q)$. The
$(-p)$-th vertical line of $\Gamma^{\bullet}(C^{\bullet})$ is
$\Gamma^{-p}(E^{\bullet})$ and $\Gamma^{-p}(-)=\oplus_{
i_0<\dots<i_p}\Gamma_{I_{i_0}+\dots+I_{i_p}}(-)$. Hence the vertical
cohomology of
$\Gamma^{\bullet}(C^{\bullet})$ in degree $(-p,q)$ is $\oplus_{
i_0<\dots<i_p}H^q(\Gamma_{I_{i_0}+\dots+I_{i_p}}(E^{\bullet}))$.
 But
$H^q(\Gamma_{I_{i_0}+\dots+I_{i_p}}(E^{\bullet}))\cong
H^q_{I_{i_0}+\dots+ I_{i_p}}(M)$. Hence in the corresponding
spectral sequence
$E_1^{-p,q}=\oplus_{i_0<\dots<i_p}H^q_{I_{i_0}+\dots+ I_{i_p}}(M)$
and the total degree of $E_1^{-p,q}$ is $t=q-p$.
\qed

\begin{remark}\label{E_2}
The differentials on the $E_1$ page have bidegree
$(1,0)$, hence $E_2^{-p,q}$ is the
$(-p)$-th cohomology of the resulting complex
$$0\to E_1^{-n+1,q}\to E_1^{-n+2,q}\to \dots\to E_1^{-1,q}\to
E_1^{0,q}\to 0$$
where the differential $E_1^{-p,q}\to E_1^{-p+1,q}$ takes $x\in
H^q_{I_{i_0}+\dots+ I_{i_p}}(M)\subset E_1^{-p,q}$ to
$\oplus_{j=0}^{j=p}(-1)^jh^q_{J, J_j}(x)$. Here
$J_j=I_{i_0}+\dots+I_{i_{j-1}}+I_{i_{j+1}}+\dots+I_{i_p}$ and
$J=I_{i_0}+\dots+I_{i_p}$ while $h^q_{J,J_j}:H^q_{J}(M)\to
H^q_{J_j}(M)$ is the natural map induced by the inclusion
$J_j\subset J$ (hence $h^q_{J, J_j}(x)\in H^q_{J_j}(M)\subset
E_1^{-p+1,q}$).
\end{remark}

It is not hard to see that if $n=2$, i.e. there are just two ideals, $I_1$
and
$I_2$, the Mayer-Vietoris spectral sequence becomes the standard
Mayer-Vietoris long exact sequence. Furthermore, it is not hard to see,
that if $n\leq 3$, i.e. there are at most three ideals, then the
Mayer-Vietoris spectral sequence degenerates at $E_2$. 

\medskip

\noindent {\bf Open question.} {\it Does the Mayer-Vietoris spectral
sequence always degenerate at
$E_2$? Does it degenerate at $E_2$ at least when $A$ is regular and
$M=A$?}

\medskip

It is shown in \cite[1.2ii]{AGZ} that the answer is positive if $A$ and
$A/I_i$ are all regular and $M=A$, but the proof does not extend to more
general ideals $I_i$, even if $A$ is still regular.

\medskip

The following two corollaries, \ref{plan} and \ref{MVDelta}, play a key
role in the proofs of our main results. They show that under suitable
assumptions on integers $d$ and
$t$ the Mayer-Vietoris spectral sequence degenerates at
$E_2$ in degree $(-t,d)$, the abutment $H_I^{d-t}(M)$ is isomorphic 
to $E_2^{-t,d}$ and 
is expressible in terms of the reduced singular homology of a
certain simplicial complex. For integers $d$ and $t$ from the statement
of Theorem \ref{maintheorem} the above-mentioned assumptions will be
proven in Theorem \ref{main} in the next section.

\begin{corollary}\label{plan}
Let $A$ be a commutative ring, let $I_1,\dots,I_n\subset A$ be ideals,
let $I=I_1\cap\dots\cap I_n$, let $M$ be an $A$-module with ${\rm
dimSupp}M\leq d$ and let 
$$E_1^{-p,q}=\oplus_{i_0<\dots<i_p}H^q_{I_{i_0}+\dots+I_{i_p}}(M)
\Longrightarrow H^{q-p}_{I}(M)$$
be the associated Mayer-Vietoris spectral sequence. Let $t<d$ be a
non-negative integer. Assume that if $p<t$, then 
$H^q_{I_{i_0}+\dots+I_{i_p}}(M)=0$ for all  
$q\geq d-t+p$ and all $\{i_0,\dots,i_p\}$. Then the 
spectral sequence degenerates at $E_2$ in degree $(-t,d)$, i.e.
$E_2^{-t,d}=E_{\infty}^{-t,d}$ and
$H^{d-t}_I(M)\cong E_2^{-t,d}$.
\end{corollary}

\emph{Proof.} The terms
$E_1^{-p,q}=\oplus_{j_0<\dots<j_p}H^q_{I_{j_0}+\dots+I_{j_p}}(M)$ 
are zero for $q>d$ since $H^q_J(M)=0$ for any ideal
$J$ provided $q>d\geq{\rm dimSupp}M$.

The terms $E_1^{-p,q}$ with
$q-p=d-t$ are zero if $0\leq p<t$ for then  
$H^{q=d-t+p}_{I_{j_0}+\dots+I_{j_p}}(M)=0$ by assumption. They are also
zero if $p>t$ for then $q=d-t+p>d$. Thus the only possibly non-zero 
$E_1^{-p,q}$ with $q-p=d-t$ is
$E_1^{-t,d}$. Hence
$H^{d-t}_I(M)\cong E_{\infty}^{-t,d}$.

If $r\geq 2$, the incoming differentials $d^r:E_r^{-t-r,d+r-1}\to
E_r^{-t,d}$ are zero because $E_r^{-t-r,d+r-1}$ is a subquotient of
$E_1^{-t-r,d+r-1}$ which vanishes since $d+r-1>d$. 

If $r\geq 1$, the outgoing differentials
$d^r:E^{-t,d}_r\to E^{-t+r,d-(r-1)}_r$ are also zero because
$E_1^{-t+r,d-(r-1)}=0$ by assumption (in this case $p=t-r<t$ and
$q=d-(r-1)$, hence $q-p=d-t+1>d-t$).

Hence $E_{\infty}^{-t,d}\cong E_2^{-t,d}$.\qed

\begin{remark}\label{(-t,d)}
As is shown above, the differential $d^1:E^{-t,d}_1\to
E^{-t+1,d}_1$ is zero, hence $E_2^{-t,d}$ is the cokernel of the map
$E_1^{-t-1,d}\to E_1^{-t,d}$ of the complex $E_1^{-\bullet,d}$ of Remark
\ref{E_2}.
\end{remark}

\begin{corollary}\label{MVDelta}
Let $A$ be a complete local ring with maximal ideal $\mathfrak m$. Let
$I=I_1\cap I_2\cap\dots\cap I_n$ be the intersection of $n$ ideals 
$I_1, \dots,I_n$ of $A$. Let $P$ be a prime ideal of $A$ such that ${\rm
dim}(A/P)=d$ and let $M$ be an
$A$-module supported on $V(P)$. Let 
$$E_1^{-p,q}=\oplus_{i_0<\dots<i_p}H^q_{I_{i_0}+\dots+I_{i_p}}(M)
\Longrightarrow H^{q-p}_{I_1\cap\dots\cap I_n}(M)$$
be the associated Mayer-Vietoris spectral sequence. Let
$\Delta$ be the simplicial complex on $n$ vertices $\{1,2,\dots,n\}$
defined as follows: a simplex
$\{j_0,\dots,j_p\}$ belongs to $\Delta$ iff
$I_{j_0}+\dots+I_{j_p}+P$ is not $\mathfrak m$-primary. For any integer
$t\geq 0$ the $E_2^{-t,d}$ term is isomorphic to 
$\tilde H_{t-1}(\Delta;H^d_{\mathfrak m}(M))$, the
$(t-1)$-th reduced singular homology group of $\Delta$ with
coefficients in
$H^d_{\mathfrak m}(M)$.
\end{corollary}

\emph{Proof.} By the Hartshorne-Lichtenbaum local vanishing theorem
\cite[3.1]{Ha}, 
$H^d_J(M)\ne 0$ only if
$\sqrt{J+P}=\mathfrak m$ in which case $H^d_J(M)\cong H^d_{\mathfrak
m}(M)$. Thus 
$H^d_{I_{j_0}+\dots+I_{j_p}}(M)\ne 0$ only if $\{j_0,\dots,j_p\}$ is
{\it not} a simplex of $\Delta$, in which case 
$H^d_{I_{j_0}+\dots+I_{j_p}}(M)\cong H^d_{\mathfrak m}(M)$. The natural
map
$h^d_{J,J_j}:H^d_{J}(M)\to
H^d_{J_j}(M)$ of Remark \ref{E_2} is non-zero
only if $\{i_0,\dots,\widehat{i_j},\dots,i_{p}\}$ is {\it not} a
simplex of 
$\Delta$ in which case it is the identity map. 
Thus the (cohomological) complex 
$E_1^{-\bullet,d}$ of Remark \ref{E_2}, upon giving $E_1^{-p,d}$ the
name $C_p$, turns into a (homological) complex $C_{\bullet}$ which is the 
standard complex for the
computation of the singular homology of the pair
$(S,\Delta)$ with coefficients in the module $H^d_{\mathfrak m}(M)$,
where $S$ is the full
$(n-1)$-simplex on the set $\{1,\dots,n\}$. Thus 
$E_2^{-t,d}=H^{-t}(E_1^{-\bullet,d})=H_t(C_{\bullet})=H^{\rm
sing}_t(S,\Delta;H^d_{\mathfrak m}(M))$. There is a long exact sequence
$\dots\to H_{s+1}(S,\Delta;G)\to H_s(\Delta;G)\to H_s(S;G)\to
H_{s}(S,\Delta;G)\to\dots$ for every abelian group $G$ which implies (upon
setting $G=H^d_{\mathfrak m}(M))$, that $H^{\rm sing}_t(S,\Delta;
H^d_{\mathfrak m}(M))\cong
\tilde H_{t-1}(\Delta;H^d_{\mathfrak m}(M))$ since $S$ is
contractible.\qed

\medskip

The plan of our proof of Theorem \ref{maintheorem} and the more general
results of Section 4 is as follows. In the next section we will
prove Theorem
\ref{main} which among other things implies that if $d$ and $t$ are
integers from the statement of Theorem \ref{maintheorem} and $p<t$, then
$H^q_{I_{i_0}+\dots+I_{i_p}}(M)=0$ for all  
$q\geq d-t+p$ and all $\{i_0,\dots,i_p\}$. Corollary \ref{plan} then
implies that
$H^{v+1}_I(M)\cong E_2^{-t,d}$, where $v=d-t-1$, as in the statement of
Theorem
\ref{maintheorem}. Now Corollary \ref{MVDelta} (with $P=0$) yields
$E_2^{-t,d}=\tilde H_{t-1}^{\rm
sing}(\Delta;H^d_{\mathfrak m}(M))\cong H^d_{\mathfrak m}(M)\otimes_k
\tilde H_{t-1}^{\rm
sing}(\Delta;k)$, where $k$ is a coefficient field of $R$ (the last
isomorphism holds by the universal coefficients theorem, since we are over
a field). This completes the proof of Theorem \ref{maintheorem} modulo 
Theorem
\ref{main}.

But in the special case that $t=[\frac{d-2}{c}]= 2$ (the smallest value
of $t$ for which Theorem \ref{maintheorem} was not previously known) we
can complete the proof of Theorem
\ref{maintheorem} very quickly without appealing to 
Theorem \ref{main}. Namely, result
(ii) of the Introduction shows that
$H^{q}_{I_{j}}(M)=0$ for $q\geq d-2$, i.e. the hypotheses of Corollary
\ref{plan} hold for $p=0$. In the next two paragraphs we will prove that
$H^{q}_{I_{i_0}+I_{i_1}}(M)=0$ for $q\geq d-1$, i.e. 
the hypotheses of Corollary \ref{plan} hold for $p=1$ as well. Thus the
hypotheses of Corollary \ref{plan} hold for $p<t=2$ and the above proof
of Theorem \ref{maintheorem} goes through without any reference to
Theorem \ref{main}.

It remains to show that $H^{q}_{I_{i_0}+I_{i_1}}(M)=0$ for $q\geq d-1$.
Equivalently, it remains to show that every minimal prime of
$I_{i_0}+I_{i_1}$ has dimension at least two and the punctured
spectrum of $R/(I_{i_0}+I_{i_1})$ is connected (this equivalence is
due to Ogus \cite[2.11]{O} in characteristic zero and Peskine and Szpiro
\cite[III, 5.5]{PSz} in characteristic $p>0$; a characteristic-free proof
has been given by Huneke and Lyubeznik \cite[2.9]{HL}). 

Since $R$ is regular, the height of every
minimal prime over the ideal $I_{i_0}+I_{i_1}$ is at most the sum of the
heights of $I_{i_0}$ and $I_{i_1}$, ie. $2c\leq d-2$. Hence the dimension
of every minimal prime over $I_{i_0}+I_{i_1}$ is at least 2. As for the
connectedness, let
$k$ be a coefficient field of $R$ and let $X_1,\dots, X_d$ generate
$\mathfrak m$. Then
$R/(I_{i_0}+I_{i_1})\cong ((R/I_{i_0})\hat \otimes_k(R/I_{i_1}))/{\rm
Diag}$, where $\hat \otimes_k$ is the complete tensor product over $k$
and Diag is the ideal generated by the $d$ elements $X_i\otimes 1 -
1\otimes X_i$. The ring
$(R/I_{i_0})\hat \otimes_k(R/I_{i_1})$ has a unique minimal prime since
$k$ is separably closed \cite[4.5]{HL}. The sum of the heights  of
$I_{i_0}$ and
$I_{i_1}$ is at most $2c\leq d-2$, so the dimension of $(R/I_{i_0})\hat
\otimes_k(R/I_{i_1})$ is at least $d+2$, hence the Faltings connectedness
theorem \cite{Fa1}, \cite{Fa2}, \cite[3.1, 3.3]{HoHu} implies that the
punctured spectrum of
$R/(I_{i_0}+I_{i_1})\cong ((R/I_{i_0})\hat \otimes_k(R/I_{i_1}))/{\rm
Diag}$ is connected.   \qed

\medskip

We conclude this section by showing that the case
$t=[\frac{d-2}{c}]= 2$ of Theorem \ref{maintheorem} that has just been
proven quickly leads to a complete understanding of an example that proved
intractable in our old paper \cite[5.12]{HL} by means of techniques that
were known to us then. In that old paper we, assuming that
$[(d-2)/c]=2$, completely analyzed the cases where $I$ has at most
five minimal primes by representing $I$ as the intersection of two
ideals and then using the standard Mayer-Vietoris long exact sequence
\cite[5.11]{HL}. But for six minimal primes we found the following example
for which the Mayer-Vietoris long exact sequence failed to provide an
answer. Below we show that the Mayer-Vietoris spectral sequence 
 works beautifully for
this example (via
the
$[(d-2)/c]= 2$ case of Theorem \ref{maintheorem}) and produces an
unexpected result: cd$(R,I)$ depends on the characteristic of the field
$k$.

\begin{example}\label{exa}
Let $R$ be a complete $d$-dimensional equicharacteristic regular local
 ring with maximal ideal $\mathfrak m$ and a separably closed residue
field $k$. Let $I$ be an ideal of $R$
with six minimal primes $I_1,\dots I_6$ such that the height of each
$I_i$ equals $c$ where $2c+2\leq d$. Set $\Lambda=\{\{1,2,3\}, \{1,3,4\},
\{2,4,5\},\{2,4,6\}$, $\{1,5,6\},$ $\{3,5,6\}, \{2,3,5\}, \{3,4,6\},
\{1,2,6\}, \{1,4,5\}\}.$ Suppose $\{i,j,\kappa\}\in \Lambda$ if and only
if
$I_i+I_j+I_{\kappa}$ is $\mathfrak m$-primary. Then $H^{d-2}_I(R)=0$ if
and only if ${\rm char}k\ne 2$ and
${\rm cd}(R,I)\leq d-3$ also if and only if ${\rm char}k\ne 2$.
\end{example}

\emph{Proof.} It is not hard to check that every four-element subset of
$\{1,\dots,6\}$ contains some element of $\Lambda$ as a subset. This
means that the complex $\Delta$ of Theorem \ref{maintheorem} does not
contain any 3-simplicies, i.e. it is a union of 2-simplices. A 2-simplex
$\{i,j,\kappa\}$ belongs to $\Delta$ if and only if it is not an element
of
$\Lambda$. Hence the complete list of the 2-simplices of $\Delta$ is 
$\{\{1,2,4\}, \{1,2,5\}$, $\{1,3,5\}, \{1,3,6\}, \{1,4,6\},
\{2,3,4\}, \{2,3,6\}, \{3,4,5\}, \{4,5,6\}, \{2,5,6\}\}$ and the picture
below shows that this collection of simplicies triangulates the real
projective plane (which is obtained from the 2-dimensional disc (or
 hexagon) by identifying pairs of antipodal points on the
boundary).

\medskip

\begin{picture}(300,100)
\put(150,100){\line(1,-1){60}}
\put(90,40){\line(1,1){60}}
\put(151,100){5}
\put(150,100){\circle*{4}}
\put(210,40){3}
\put(210,40){\line(-2,-1){60}}
\put(150,2){5}
\put(150,10){\circle*{4}}
\put(210,40){\line(0,1){30}}
\put(210,70){\line(-2,1){60}}
\put(210,40){\circle*{4}}
\put(84, 40){6}
\put(90,40){\circle*{4}}
\put(90,40){\line(1,0){120}}
\put(90,40){\line(0,1){30}}
\put(90,70){\line(2,1){60}}
\put(117,72){1}
\put(90,70){\circle*{4}}
\put(84,70){3}
\put(120,70){\circle*{4}}
\put(180, 70){2}
\put(210,70){\circle*{4}}
\put(210,70){6}
\put(90,40){\line(2,-1){60}}
\put(180,70){\circle*{4}}
\put(147,42){4}
\put(150,40){\circle*{4}}
\put(120,70){\line(1,0){60}}
\put(120,70){\line(1,-1){30}}
\put(150,40){\line(1,1){30}}
\put(150,10){\line(0,1){30}}
\put(90,70){\line(1,0){30}}
\put(180,70){\line(1,0){30}}
\end{picture}

\medskip

Hence the complex
$\Delta$ is homeomorphic to the real projective plane. Thus 
$\tilde H_1(\Delta;k)$, the critical singular homology group of
$\Delta$, is zero if the charateristic of $k$ is different from 2 and is
isomorphic to $k$ otherwise. Now we are done by the $[(d-2)/c]=2$ case of
Theorem
\ref{maintheorem}.\qed

\medskip

For a concrete realization of the above example let $a\in k$ be an
element such that $a\ne 0,1,-1$ and $a^2+a-1\ne 0$. Let
$R=k[[X_1,X_2,X_3,X_4,X_5,X_6]]$ be the ring of formal power series in
six variables $X_1,\dots, X_6$. Let ideals $I_1\dots, I_6$ of height two
be defined as follows: 

$I_1=(X_1,\ X_2),$ 
$I_2=(X_3, \ X_4),$ 
$I_3= (X_5,
\ X_6),$ 

$I_4= (X_1+X_3+X_6, \ \ \frac{1}{a^2+a-1}X_1+X_4+\frac{1}{a}X_6),$

$I_5=(X_1+X_3+X_5, \ \ X_2+aX_3+X_5),$ 

$I_6=(X_2+X_4+X_5,\ \ 
\frac{1}{a}X_2+aX_4+X_6)$. 

It is tedious but straightforward to verify
that $I_i+I_j+I_{\kappa}$ is $\mathfrak m$-primary iff $\{i,j,\kappa\}\in
\Lambda$ where $\Lambda$ is the same as in the above example (the
assumptions $a\ne 1, -1$, and their consequences
$\frac{1}{a}\ne a, \frac{1}{a^2+a-1}$ are used in the
verification). Setting
\hbox{$I=I_1\cap\dots\cap I_6$} we conclude that $H^4_I(R)=0$ if
char$k\ne 2$ and $H^4_I(R)\ne 0$ if char$k=2$. Also, cd$(R,I)\leq 3$ 
if char$k\ne 2$ and cd$(R,I)=4$ if char $k=2$.

\section{A vanishing theorem for $H^q_{I_0+\dots+I_{\mathfrak p}}(M)$.}
\label{S:QC} The main result of this section is Theorem \ref{main} which
establishes the vanishing of $E_1^{-p,q}$ needed in the
proofs of our main results. We prove our main results
for an arbitrary commutative Noetherian local ring containing a field.
 This necessitates working not with the
height of the minimal primes over $I$, like in the statement of
Theorem \ref{maintheorem}, but with the integer $c(I)$
whose definition we recall below and which is
equal to the maximum height of a minimal prime over $I$ if the
ring is regular.

\medskip

\noindent {\bf Definition.} \cite[2.1]{HL} {\it If
$A$ is a local ring and
$I\subset A$ is an ideal, 
$$c(I)={\rm emb.dim}A-{\rm min}\{{\rm dim}(A/P)|P\ \text {is a minimal
prime over }I\}.$$}
In the sequel
$\hat A$ denotes the completion of
$A$ with respect to the maximal ideal. The following result is a
rephrasing of Faltings
\cite[Korollar 2]{Fa}.
\begin{theorem}\label{Falt}
Let $A$ be a local ring containing a field, let $c>0$ and
$d\geq 0$ be integers, let $I\subset A$ be an ideal with
$c(I\hat A)\leq c$  and let $M$ be an $A$-module such
that ${\rm dimSupp}M\leq d$. Then
$H^q_I(M)=0$ for $q>d-[\frac{d-1}{c}]$.
\end{theorem}
\emph{Proof.} 
\cite[Korollar 2]{Fa} says that if $M$ is finitely generated,
then
$H^q_I(M)=0$ for $q>(1-\frac{1}{c(I\hat A)})\cdot{\rm dim}M+1$ provided
$c(I\hat A)\geq 2$ and $H^q_I(M)=0$ for \hbox{$q>c(I\hat A)$} provided
$c(I\hat A)<2$. This immediately implies Theorem \ref{Falt}
because $M$ is the direct limit of its finitely generated submodules, the
dimension of each finitely generated submodule is at most $d$ and 
local cohomology commutes with direct
limts, and for all integers
$d\geq 0$ and $c>0$ the function
$\phi(c,d)=d-[\frac{d-1}{c}]=[(1-1/c)d+1]$ is non-decreasing
and 
$\phi(c,d)\geq 1$. \qed

\medskip
 
We need the following generalization.
\begin{theorem}\label{arbitr}
Let $A$ be a local ring containing a field, let $c$ and $d$ be positive
integers, let
$J=I_0+\dots+ I_{\mathfrak p}$ be the sum of $\mathfrak p+1$ ideals $I_0,\dots, I_{\mathfrak p}$
such that
$c(I_j\hat A)\leq c$ for each $j$ and let $M$ be an $A$-module
such that ${\rm dimSupp}M\leq d$. Then 
$H^q_J(M)=0$ for
$q>d-[(d-1)/c]+\mathfrak p$. In particular, if dim$A\leq d$,
then $cd(A,I)\leq d-[(d-1)/c]+\mathfrak p$.
\end{theorem}
\emph{Proof.} We use induction on $\mathfrak p$, the case $\mathfrak p=0$ being known by
Theorem \ref{Falt}. 
Assume $\mathfrak p>0$ and the result proven for $\mathfrak p-1$. Since a prime contains 
the ideal $I_{\mathfrak p}\cap(I_0+\dots+I_{\mathfrak p -1})$ if and only if it contains 
$\mathcal J=I_{\mathfrak p}\cap I_0+I_{\mathfrak p}\cap I_1+\dots+I_{\mathfrak p}\cap I_{\mathfrak p -1}$,
these two ideals have the same radical, so there is an isomorphism of the
corresponding local cohomology functors
$H^*_{I_{\mathfrak p}\cap(I_0+\dots+I_{\mathfrak p -1})}(-)\cong H^*_{\mathcal J}(-)$. The ideal
$\mathcal J$ (resp. $I_0+\dots+I_{\mathfrak p -1}$) is the sum of $\mathfrak p$ ideals $I_{\mathfrak p}\cap
I_0,\dots, I_{\mathfrak p}\cap I_{\mathfrak p -1}$ (resp. $I_0,\dots,I_{\mathfrak p -1}$) 
such that $c((I_{\mathfrak p}\cap I_j)\hat A)\leq c$ (resp. 
$c(I_j)\leq c$) for every $j\leq \mathfrak p -1$. Hence by the induction
hypothesis $H^q_{\mathcal J}(M)\cong
H^q_{I_{\mathfrak p}\cap(I_0+\dots+I_{\mathfrak p -1})}(M)\cong 0$ and
$H^q_{I_0+\dots+I_{\mathfrak p -1}}(M)=0$ for 
$q>d-[(d-1)/c]+\mathfrak p -1$. The 
Mayer-Vietoris exact sequence 
$H^{q-1}_{I_{\mathfrak p}\cap(I_0+\dots+I_{\mathfrak p -1})}
(M)\to H^q_J(M)\to
H^q_{I_{\mathfrak p}}(M)\oplus H^q_{I_0+\dots+I_{\mathfrak p -1}}(M)$ 
implies the
theorem.\qed

\medskip

The following lemma is a consequence of \cite[2.2, 2.3]{HL}:
\begin{lemma}\label{c(I)}
Let $A$ be a complete local ring containing a field and let
$I$ be an ideal of $A$. Then

(i) $c(I)=c(IB)$ where
$B$ is the completion of the strict Henselization of $A$ (i.e. 
$B\cong K\hat\otimes_k\hat A$, where $k$ is a coefficient field of
$A$, while
$\hat\otimes_k$ is the complete tensor product over $k$ and $K$ is the
separable closure of $k$).

(ii) $c(I\widehat{A_P})\leq c(I)$ for every prime ideal $P\supset I$,
where $\widehat{A_P}$ is the completion of $A_P$.
\end{lemma}
\emph{Proof.} \cite[2.2]{HL} says that if $A$ is a universally catenary
local ring containing a field, then $c(I)=c(I\hat A)$ and $c(I)=c(IB)$
where $B$ is the completion of the strict Henselization of $\hat A$ while
\cite[2.3]{HL} says that under the same assumptions $c(I_P)\leq c(I)$.
But a complete local ring $A$ is universally catenary and so are all
its localizations $A_P$.
\qed

\medskip 

The following result is a rephrasing of 
Huneke and Lyubeznik
\cite[3.8]{HL}.
\begin{theorem}\label{HL}
Let $A$ be a local ring containing a field, let $B=\widehat{\hat A^{\rm
sh}}$ be the completion of the strict Henselization of the completion
of $A$, let
$c$ and
$d$ be positive integers with $c<d$, let $I\subset A$ be an ideal such
that $c(I\hat A)\leq c$ and $\sqrt{IB}$ is a prime ideal of $B$ and
 let
$M$ be an
$A$-module with ${\rm dimSupp}M\leq d$. Then $H^q_I(M)=0$ for
$q>d-1-[\frac{d-2}{c}]$.
\end{theorem}
The actual statement of \cite[3.8]{HL} is a
special case with $M$ finitely generated, $0<c(I\hat A)=c$ and $d={\rm
dim}M$. Theorem
\ref{HL} is immediate from this special case just like Theorem
\ref{Falt} is immediate from \cite[Korollar 2]{Fa}.
The following generalization is the main result of this section.
\begin{theorem} \label{main}
Let $A$ be a local ring containing a field, let $B=\widehat{\hat A^{\rm
sh}}$ be the completion of the strict Henselization of the completion
of $A$, let $c$ be a
positive integer and let 
$J$ be an ideal of $A$ such that 
$\sqrt{JB}=\sqrt{I_0+\dots+ I_{\mathfrak p}}$ where
 $I_0,\dots,
I_{\mathfrak p}$ are prime ideals of $B$ with $c(I_j)\leq c$ 
for each $j$.
Let
$M$ be an
$A$-module with ${\rm dimSupp}M\leq d$  where $d>(\mathfrak p +1)c$. Then 
$H^q_J(M)=0$ for $q>d-1-[(d-2)/c]+\mathfrak p .$ 
In particular, if $(\mathfrak p +1)c<{\rm dim}A$ and ${\rm dim}A\leq d$,
then
$cd(A,I)\leq d-1-[(d-2)/c]+\mathfrak p $.
\end{theorem}

{\it Remarks.} (i) The restriction $d>c(\mathfrak p +1)$ cannot be omitted. For
example, if
$A$ is regular of dimension $d=c(\mathfrak p +1)$ and $I_0,\dots, I_{\mathfrak p}$ each have
height
$c$ and add up to the maximal ideal, then $H^d_J(A)\ne 0$ while
$d>d-1-[(d-2)/c]+\mathfrak p$.

(ii) The condition $\mathfrak p<t=[\frac{d-2}{c}]$ implies
$d>c(\mathfrak p +1)$, hence Theorem \ref{main} does indeed show that the
hypotheses of Corollary
\ref{plan} hold for $t=[\frac{d-2}{c}]$.

\medskip

Our proof of Theorem \ref{main} is considerably longer than the above
proof of Theorem \ref{arbitr}. One cannot deduce Theorem \ref{main}
from Theorem \ref{HL} in the same way as Theorem \ref{arbitr} was deduced
from Theorem \ref{Falt} because the ideals $I_{\mathfrak p}\cap
I_j$ are not formally geometrically irreducible, hence the induction
hypotheses does not apply to $\mathcal J=I_{\mathfrak p}\cap I_0+\dots+I_{\mathfrak p}\cap
I_{\mathfrak p -1}$.

\medskip

\emph{Proof of Theorem \ref{main}.}
The module $H^q_{JB}(B\otimes_AM)\cong B\otimes_AH^q_J(M)$ vanishes if and
only if
$H^q_J(M)$ vanishes since $B$ is faithfully flat over $A$. 
Since \hbox{dimSupp$B\otimes_AM\leq d$,} replacing the ring
$A$ and the module $M$ by $B$ and 
$B\otimes_AM$ respectively, we may assume that $A$ is complete with
separably closed residue field and the $I_i$s are prime ideals of $B$.

By one of Cohen's structure theorems for complete local rings,
the ring $A$ contains a coefficient field $k\subset A$
and there exists a surjective $k$-algebra homomorphism $R\to A$ where
$R=k[[X_1,\dots,X_n]]$ is the ring of formal power series in $n$
variables over $k$ and $n$ is the embedding dimension of $A$. The
$A$-module
$M$ acquires a structure of
$R$-module via this homomorphism. If
$I\subset A$ is an ideal and
$\tilde I\subset R$ is the full preimage of $I$, then \hbox{$H^i_I(M)\cong
H^i_{\tilde I}(M)$}. Each of the ideals $\tilde I_i$ is prime and
$c(\tilde I_i)=c(I_i)$. Thus replacing $A$ by $R$ and $I_i$ by $\tilde
I_i$ we may assume that $A=k[[X_1,\dots,X_n]]$ is a ring of formal power
series over a separably closed field $k$. 

Let $K$ be
an uncountable separably closed field extension of $k$. The 
inclusion
$A=k[[X_1,\dots,X_n]]\hookrightarrow K[[X_1,\dots,X_n]]=R'$ induced
by the natural inclusion $k\hookrightarrow K$ makes
$R'$ an $A$-algebra. Using the Koszul resolution of $A/m\cong k$ on
$X_1,\dots, X_n$ one sees that Tor$_1^{A}(k, R')\cong 0$, so $R'$
is flat over $A$ by the local criterion of flatness
\cite[20.C(3$'$)]{Mats}. Since the inclusion $A\hookrightarrow
R'$ is a local ring homomorphism, $R'$
is faithfully flat over $A$ \cite[4.D]{Mats}. Clearly, $c(I_iR')=c(I_i)$
and dimSupp$(R'\otimes_AM)={\rm dimSupp}M\leq d$. At this point we recall
the statement of \cite[4.3]{HL} which in our notation says the following:
{\it if
$k$ is separably closed in $K$, then for every prime ideal $P\subset A$
the radical of the ideal $PR'$ of $R'$ is prime.} Since in our case $k$ is
a separably closed field, the radicals of the ideals
$I_iR'$ are prime. Since
$R'$ is flat \hbox{over $A$,} $H^i_{IR'}(R'\otimes_AM)\cong
R'\otimes_AH^i_I(M)$ for every ideal $I$ of $A$ and since $R'$ is
faithfully flat,
$H^i_{IR'}(R'\otimes_AM)\cong 0$ if and only if $H^i_I(M)\cong 0$. Hence 
replacing the ring $A$, the ideal $I_i$ and the module
$M$ by $R'$, rad($I_iR'$) and
$R'\otimes_AM$ respectively, we may assume that $A=K[[X_1,\dots,X_n]]$ is
a ring of formal power series over an uncountable separably closed field
$K$ and each $I_i$ is prime.

We may assume
that $M$ is finitely generated because $H^q_J(M)$ is the direct limit of
$H^q_J(M')$ as $M'$ runs over all the finitely generated submodules of
$M$ and dim$M'\leq d$ for all such $M'$. 

Let $v=d-1-[(d-2)/c]+\mathfrak p$. We may assume that $c\not |(d-1)$ because
otherwise
$v=d-[(d-1)/c]+\mathfrak p$ and we are done by Theorem \ref{arbitr}. 

We use induction on $\mathfrak p$, the case
$\mathfrak p=0$ being known by Theorem \ref{HL}. Assume $\mathfrak p>0$ and let
$J'=I_0+\dots+I_{\mathfrak p -1}$. The composition of functors
$\Gamma_{I_{\mathfrak p}}(\Gamma_{J'}(-))\cong \Gamma_J(-)$ leads to the spectral
sequence 
$$E^{p,q}_2=H^p_{I_{\mathfrak p}}(H^q_{J'}(M))\Longrightarrow H^{p+q}_J(M).$$
To prove the theorem it is enough to show that $E^{p,q}_{\infty}=0$
provided $p+q>v$, and we are going to show this. For some, but not all,
pairs
$(p,q)$ with
$p+q>v$ we are even going to show that
$E_2^{p,q}=0$.

To this end, let $p$ and $q$ be non-negative integers with
$p+q>v.$ There are only two possibilities, either
${\rm dimSupp}H^q_{J'}(M)>c$, or ${\rm dimSupp}H^q_{J'}(M)\leq c$.

First we consider the case that ${\rm dimSupp}H^q_{J'}(M)>c$. We are
going to prove that $E_2^{p,q}=0$ in this case. Let
$\delta$ be the biggest integer such that $q\leq
d-\delta-[(d-\delta-1)/c]+\mathfrak p -1$ (it exists since
$d-x-[(d-x-1)/c]+\mathfrak p -1$ is a non-increasing function of $x$). If 
$\delta'>\delta$, then
$q>d-\delta'-[(d-\delta'-1)/c]+\mathfrak p -1$ and if $P$ is a
prime of
$A$ of dimension $\delta'$, then $c(I\widehat{A_P})\leq c$ by Lemma
\ref{c(I)} and dimSupp$M_P\leq d-\delta'$, so
$H^q_{J'A_P}(M_P)=0$ by Theorem \ref{arbitr}. Hence the
module $H^q_{J'}(M)$ vanishes at all primes of dimension bigger than
$\delta$, which implies that 
$\delta\geq{\rm dimSupp}H^q_{J'}(M)>c$. Now Theorem \ref{HL}
implies that $E^{p,q}_2=H^p_{I_{\mathfrak p}}(H^q_{J'}(M))=0$ for
$p>\delta-1-[(\delta-2)/c]$. Hence it
is enough to show that $p>\delta-1-[(\delta-2)/c]$.  

This is a formal consequence of the inequalities 
$p+q>d-1-[(d-2)/c]+\mathfrak p$ and
$q\leq d-\delta-[(d-\delta-1)/c]+\mathfrak p -1$, which imply
\begin{align*}
p& >(d-1-[(d-2)/c]+\mathfrak p )-q \\
& \geq (d-1-[(d-2)/c]+\mathfrak p )-(d-\delta-[(d-\delta-1)/c]+\mathfrak p -1)\\
&=\delta+[(d-\delta-1)/c]-[(d-2)/c],
\end{align*}
 so it is enough to show that
$$\delta+[(d-\delta-1)/c]-[(d-2)/c]\geq \delta-1-[(\delta-2)/c],$$ i.e.
$$[(d-\delta-1)/c]+[(\delta-2)/c]+1-[(d-2)/c]\geq 0.$$ Let
$\delta=\kappa c+r$, where $0\leq r<c$. Then
$$[(d-\delta-1)/c]+[(\delta-2)/c]=[(d-r-1)/c]+[(r-2)/c].$$ If $r\geq 2$,
then $[(r-2)/c]=0$, so we need to show that $$[(d-r-1)/c]+1\geq
[(d-2)/c].$$ But $r\leq c-1$ implies $$[(d-r-1)/c]+1\geq
[(d-c)/c]+1=[d/c]\geq [(d-2)/c],$$ so we are done in this case. If
$r=0,1$, then $[(r-2)/c]=-1$, so $$[(d-r-1)/c]+[(r-2)/c]+1=[(d-r-1)/c]\geq
[(d-2)/c].$$ This completes the proof that $E_2^{p,q}=0$ provided
${\rm dimSupp}H^q_{J'}(M)>c$ and $p+q>v$.

It remains to prove that $E_{\infty}^{p,q}=0$ if 
${\rm dimSupp}H^q_{J'}(M)\leq c$ and $p+q>v$. If $p>c$, then $p>{\rm
dimSupp}H^q_{J'}(M)$, hence $E_2^{p,q}=H^p_{I_{\mathfrak p}}(H^q_{J'}(M))=0$.
Thus it remains to consider the case that $p\leq c$. 

First we assume that 
$p+q\geq v+2$. We are going to show that
$E_2^{p,q}=0$ under this assumption. For this it is
enough to show, like in the preceding paragraph, that
dimSupp$H^q_{J'}(M)<p$, i.e. that
$H^q_{J'A_P}(M_P)=0$ for every prime $P$ of dimension $\geq p$. Lemma
\ref{c(I)} implies that $c(J'\widehat{A_P})\leq c$. We have that 
$v=d-[(d-1)/c]+(\mathfrak p -1)$ since $c\not | (d-1)$. Since dimSupp$M_P\leq d-p$,
Theorem \ref{arbitr} implies $H^q_{J'A_P}(M_P)=0$ for
$q>d-p-[(d-p-1)/c]+(\mathfrak p -1)$
that is (since $p\leq c$ implies that $[(d-p-1)/c]\geq [(d-1)/c]-1$) for
$q>d-[(d-1)/c]+(\mathfrak p -1)-(p-1)=v-p+1$, i.e. for $p+q>v+1$. This
concludes the proof that $E_2^{p,q}=0$ provided 
$p+q\geq v+2$.

It remains to consider the case that ${\rm dimSupp}H^q_{J'}(M)\leq
c$ while $p\leq c$ and $p+q=v+1$. This case is the hardest. In this case
we are unable to prove that $E_2^{p,q}=0$, but we are going to prove
that $E_{\infty}^{p,q}=0$ by constructing a morphism
from our spectral sequence to a different spectral sequence $\bar E$ 
about which we can prove that for the
pairs $(p,q)$ in question the induced maps $E_{\infty}^{p,q}\to \bar
E_{\infty}^{p,q}$ are isomorphisms and $\bar E_{\infty}^{p,q}=0$.

First we claim that if
$q'>v-c$, then dimSupp$H^{q'}_{J'}(M)\leq c$. Indeed, if
$P$ is a prime of dimension
$\geq c+1$, then dim$M_P\leq d-c-1$. Lemma
\ref{c(I)} implies that $c(J'\widehat{A_P})\leq c$. By Theorem
\ref{arbitr}, if 
$H^{q'}_{J'A_P}(M_P)\ne 0$, then $q'\leq
d-c-1-[(d-c-2)/c]+(\mathfrak p -1)=d-1-[(d-2)/c]+\mathfrak p -c=v-c$, i.e. $q'\leq v-c$. This
proves the claim. 

We claim there exist elements
$x_1,\dots, x_c\in I_{\mathfrak p}$ such that the ideals $P+I_{\mathfrak p}$ and
$P+(x_1,\dots,x_c)$ have the same radical for every
minimal prime $P$ of the support of $H^{q'}_{J'}(M)$ as $q'$ runs
through all integers $q'>v-c$. In fact, by induction on $c'$ we are going
to prove the following statement. {\it There exist $x_1,\dots, x_c\in I_{\mathfrak p}$
such that for each $c'\leq c$ and each minimal prime $P$ of some
$H^{q'}_{J'}(M)$, every prime $Q$ containing $P$ and $x_1,\dots,x_{c'}$
but not containing $I_{\mathfrak p}$, has dimension at most $c-c'$.} This statement
for $c'=c$ implies the claim considering that $A$ being local has only
one prime of dimension zero. 

If
$c'=0$, we set
$x_0=0$. Since dimSupp$H^{q'}_{J'}(M)\leq c$, we get dim$P\leq c$,
so the claim holds if $c'=0$. Assume $c'>0$ and
$x_1,\dots, x_{c'-1}$ have been found. To prove the claim it remains to
show that there exists an element $x_{c'}\in I_{\mathfrak p}$ which does not belong
to any prime $Q$ such that $Q$ does not contain $I_{\mathfrak p}$ and is minimal over
some ideal
$(x_1,\dots, x_{c'-1})+P$ where $P$ is minimal in the support of some
$H^{q'}_{J'}(M)$ with $q'>v-c$. For then any prime $Q'$ that does not
contain
$I_{\mathfrak p}$ and contains $(x_1,\dots, x_{c})+P$ necessarily contains some $Q$
and
$x_{c'}\not\in Q$, hence dim$Q'<{\rm dim}Q$, but dim$Q\leq c-c'+1$ by the
induction hypothesis and therefore dim$Q'\leq c-c'$, as required. The
existence of $x_{c'}$ is shown in the next two paragraphs.

The set of the minimal primes of all
the $H^{q'}_{J'}(M)$ is countable. This is because
$H^{q'}_{J'}(M)=\varinjlim{\rm Ext}^{q'}_A(A/(J')^n,M)$, so a
prime is associated to $H^{q'}_{J'}(M)$ only if it is associated to one
of the Ext$^{q'}_A(A/(J')^n,M)$, but each Ext$^{q'}_A(A/(J')^n,M)$ is
finitely generated and therefore has a finite number of associated
primes.  Thus the set of the primes $Q$ which do not contain
$I_{\mathfrak p}$ and which are minimal over some ideal $(x_1,\dots,
x_{c'-1})+P$, where $P$ is minimal in the support of some
$H^{q'}_{J'}(M)$ with $q'>v-c$, is countable, because the set of the
ideals
$P$ is countable. 

Let $y_0,\dots,y_s$ be a
finite system of generators of
$I_{\mathfrak p}$. For each element $\kappa\in K$ we set $y(\kappa)=y_0+\kappa
y_1+\kappa^2y_2+\dots+\kappa^sy_s\in I_{\mathfrak p}$. If
$\kappa_0,\kappa_1,\dots,\kappa_s\in K$ are distinct, the corresponding
Vandermond determinant is non-zero, hence $y_0,\dots,y_s$ are linear
combinations of
$y(\kappa_0),\dots,y(\kappa_s)$ with coefficients from $K$. This
implies that each prime $Q$ of our countable set of primes contains
elements
$y(\kappa)$ for at most
$s$ distinct values of $\kappa$, for otherwise
$y_0,\dots,y_s\in Q$ and hence $Q$ contains $I_{\mathfrak p}$. Thus the set of all
elements
$\kappa\in K$ such that $y(\kappa)$ belongs to some $Q$ is countable.
Since $K$ is uncountable, there exists $\kappa\in K$ such that
$y(\kappa)\in I_{\mathfrak p}$ does not belong to any $Q$. We set $x_{c'}=y(\kappa)$. 
This completes the proof of the claim. \footnote{This "countable prime
avoidance" argument is necessary because we do not know whether each
module $H^{q'}_{J'}(M)$ has only a finite number of minimal primes. If we
knew this, standard prime avoidance would do and we wouldn't have even
needed to pass from the original possibly countable field $k$ to the
uncountable field $K$. But it is still an open problem whether the set of
the minimal primes of a local cohomology module of a finitely generated
module over a Noetherian ring is always finite. One only knows that the
set of the {\it associated} primes need not be finite \cite{Kat, Si}.} 

We claim that the
Grothendieck spectral sequence of the
composition of functors $E^{p,q}_2=R^pG(R^qF(A))\Longrightarrow
R^{p+q}(G\circ F)(A)$ for fixed
$F$ and $A$ is functorial in $G$. We were unable to find a reference to
this basic fact in the literature, so we provide a proof. Let
$I^{\bullet}$ be an injective resolution of
$A$ and let $I^{\bullet,\bullet}$ be an injective double complex
resolution of $I^{\bullet}$. Let $G\to \bar G$ be a natural
transformation of functors. The spectral sequences $E$ and $\bar E$
associated to the compositions of functors $G\circ F$ and $\bar G\circ F$
respectively are the spectral sequences of the filtered total complexes
of the double complexes $G(I^{\bullet, \bullet})$ and $\bar
G(I^{\bullet,
\bullet})$ respectively. The natural transformation $G\to \bar
G$ induces a map $G(I^{p,q})\to \bar G(I^{p,q})$ for every $p$ and $q$,
hence it induces a morphism of double complexes $G(I^{\bullet,
\bullet})\to
\bar G(I^{\bullet, \bullet})$ which in turn induces a morphism
of the  total complexes that is compatible with the corresponding
filtrations. According
to
\cite[11.2.3]{G}, this gives rise to a morphism of associated spectral
sequences which proves the claim. 

We apply this claim as follows. Setting $X=(x_1,\dots,x_c)$ we get a
spectral sequence 
$$\bar E^{\pi,\kappa}_2=H^{\pi}_X(H^{\kappa}_{J'}(M))\Longrightarrow
H^{\pi+\kappa}_{X+J'}(M).$$ Since $X\subset I_{\mathfrak p}$, we have that
$\Gamma_{I_{\mathfrak p}}(N)\subset \Gamma_X(N)$ for every module $N$. This induces a
natural transformation of functors $\Gamma_{I_{\mathfrak p}}(-)\to
\Gamma_X(-)$ which according to the above claim 
induces a morphism of spectral sequences
$\phi:E\to\bar E$. 

We claim that
$\bar E^{\pi,\kappa}_2=H^{\pi}_X(H^{\kappa}_{J'}(M))\cong 0$ for
$\pi+\kappa\geq v+2$. Indeed, this is true if $\pi>c$ since $X$ is
generated by $c$ elements. If $\pi\leq c$, it has been shown above that
dimSupp$H^{\kappa}_{J'}(M)<\pi$ under the additional assumption that
$\pi+\kappa\geq v+2$. This proves the claim. 

We claim
that the resulting maps
$E_{\infty}^{p,q}\to
\bar E_{\infty}^{p,q}$ are isomorphisms if $p\leq c$ and
$p+q=v+1$. 
This claim is a special case of the following more general claim: {\it If
$p+q=v+1-s$, where $s\geq 0$ and $p+sr\leq c$, then the natural map
$E^{p,q}_r\to \bar E^{p,q}_r$ is an isomorphism.} Indeed,
if $p\leq c$ and $p+q=v+1$, then $s=0$, hence $p+sr=p\leq c$, so
$E^{p,q}_r\to \bar E^{p,q}_r$ is an isomorphism for all $r$, hence
$E^{p,q}_{\infty}\to \bar E^{p,q}_{\infty}$ is an isomorphism.

To prove this more general claim we use induction on $r$. If $r=2$,
we need to show that the natural maps  
$E^{p,q}_2=H^p_{I_{\mathfrak p}}(H^q_{J'}(M))\to H^p_{X}(H^q_{J'}(M))=\bar
E^{p,q}_2$ are isomorphisms if $p+q=v+1-s$ and $p+2s\leq c$. The
condition $p+2s\leq c$ implies that $p\leq c-2s$ and the condition
$p+q=v+1-s$ implies that $q=v+1-s-p\geq v+1-s-(c-2s)=v-c+1+s>v-c$. Hence
if $P$ is a minimal prime of $H^q_{J'}(M)$, then the ideals $P+X$ and
$P+I_{\mathfrak p}$ have the same radical. Hence if $M'$ is a finitely generated
submodule of $H^q_{J'}(M)$, then the ideals ann$M'+X$ and
ann$M'+I_{\mathfrak p}$ have the same radical, where ann$M'$ is the
annihilator of $M'$. This implies that the natural map $H^p_{I_{\mathfrak
p}}(M')\to H^p_X(M')$ is an isomorphism for every finitely generated
submodule $M'$ of $H^q_{J'}(M)$, hence the natural map
$E^{p,q}_2=H^p_{I_{\mathfrak p}}(H^q_{J'}(M))\to
H^p_{X}(H^q_{J'}(M))=\tilde E^{p,q}_2$ is an isomorphism too. This
completes the case $r=2$.

Assume the claim proven for $r-1$. Let a pair $(p,q)$ be such that
$p+q=v+1-s$ and $p+sr\leq c$ for some $s\geq 0$. We need to show that
the map $\phi_r:E_r^{p,q}\to \bar E_r^{p,q}$ is an isomorphism. We have
the following commutative diagram
$$
\CD
E^{p-(r-1),q+r-2}_{r-1}@>d_{r-1}>> E^{p,q}_{r-1} @>d_{r-1}>> E^{p+r-1,
q-r+2}_{r-1}\\ @V\phi_{r-1}VV                         
@V\phi_{r-1}VV              @V\phi_{r-1}VV\\
\bar E^{p-(r-1),q+r-2}_{r-1}@>\bar d_{r-1}>> \bar E^{p,q}_{r-1} @>\bar
d_{r-1}>>
\bar E^{p+r-1, q-r+2}_{r-1}\\
\endCD
$$
in which the vertical maps are isomorphisms by the induction hypothesis
as we are presently going to show. Indeed, for the map on the left,
setting
$p'=p-(r-1)$ and $q'=q+r-2$ we see that
$p'+q'=p+q-1=v+1-(s+1)$ and
$p'+(s+1)(r-1)=p+s(r-1)\leq p+sr\leq c$, so the induction
hypothesis holds for the map on the left. Since
$p+s(r-1)\leq p+sr\leq c$, the induction hypothesis holds
for the map in the middle. Hence the vertical maps on the left and in
the middle are isomorphisms. For the vertical map on the right we set
$p'=p+r-1$ and $q'=q-r+2$. If
$s=0$, then
$p'+ q'=p+q+1=v+2$, so both modules on the right are zero and therefore
the map on the right is an isomorphism. If $s\geq 1$, then
$p'+q'=p+q+1=v+1-(s-1)$ where $s-1\geq 0$ and 
$p'+(s-1)(r-1)=p+s(r-1)\leq p+sr\leq c$, so the induction
hypothesis holds for the map on the right. Hence the vertical map on
the right also is an isomorphism. Thus all verical maps are indeed
isomorphisms, so the middle map induces an isomorphism
$\phi_r:E^{p,q}_r\to\bar E^{p,q}_r$ on the homology modules. 
This proves the claim.

It remains to show that $\bar E^{p,q}_{\infty}=0$ if $p+q>v$.
Equivalently, it remains to show that $H^i_{X+J'}(M)=0$ for $i>v$. We
have a spectral sequence
$$\tilde E^{p',q'}_2=H^{p'}_{J'}(H^{q'}_X(M))\Longrightarrow
H^{p'+q'}_{X+J'}(M).$$
$H^{q'}_X(M)=0$ for $q'>c$ since $X$ is generated by $c$ elements. Hence
$\tilde E^{p',q'}_2=0$ if $q'>c$. 

Assume $q'\leq c$. If
$P$ is a prime of dimension $>d-q'$, then dim$M_P<q'$, so $H^{q'}_X(M)=0$.
Hence dimSupp$H^{q'}_X(M)\leq d-q'$. Since $d>c(\mathfrak p +1)$ and $q'\leq c$, we
conclude that $d-q'>c\mathfrak p$, so by induction on $\mathfrak p$,
$\tilde E^{p',q'}_2=H^{p'}_{J'}(H^{q'}_X(M))=0$ for
$p'>d-q'-1-[(d-q'-2)/c]+\mathfrak p -1$, i.e. for
$p'+q'>d-1-[(d-2)/c]+\mathfrak p -([(d-q'-2)/c]+1-[(d-2)/c])$ and hence for
$p'+q'>d-1-[(d-2)/c]+\mathfrak p =v$ as $q'\leq c$ implies
$[(d-q'-2)/c]+1-[(d-2)/c]\geq 0$.
\qed

\section{The main results}\label{S:MR}
The following theorem and its corollaries are the main results of this
paper. 

\begin{theorem}\label{general}
Let $A$ be a local ring containing a field. Let
$c>0$ and
$d>1$ be integers, let 
$I=I_1\cap\dots\cap I_n$ be the intersection of several prime ideals $I_1,
I_2,\dots$ such that $c(I_j\hat A)\leq c$ for all $j$, 
 let $t=[\frac{d-2}{c}]$, let $v=d-1-[\frac{d-2}{c}]$ and let $M$ be an
$A$-module such that ${\rm dimSupp}M\leq d$. Assume that $\sqrt
{I_jB}$ is a prime ideal of $B$ for all $j$, where $B$ is the
completion of the strict Henselization of the completion of $A$. Then 
$H^{v+1}_I(M)$ is isomorphic to the cokernel of the map
$$\Phi_{M,I}:\oplus_{j_0<\dots<j_{t +1}}H^d_{I_{j_0}+\dots+I_{j_{t
+1}}}(M)\to
\oplus_{j_0<\dots<j_{t}}H^d_{I_{j_0}+\dots+I_{j_{t}}}(M)$$
that sends every $x\in H^d_{I_{j_0}+\dots+I_{j_{t}}}(M)$ to
$\oplus_{s=0}^{s=t }(-1)^sh_s(x)$ where 
$$h_s:H^d_{I_{j_0}+\dots+I_{j_{t +1}}}(M)\to
H^d_{I_{j_0}+\dots+\widehat{I_{j_s}}+\dots+I_{j_{t +1}}}(M)$$ 
($\widehat{I_{j_s}}$ means that
$I_{j_s}$ has been omitted) is the
natural map induced by the containment $$I_{j_0}+\dots+I_{j_{t +1}}\supset
I_{j_0}+\dots+\widehat{I_{j_s}}+\dots+I_{j_{t +1}}.$$
\end{theorem}

We point out that the conditions that $d>1$ and $d$ is an integer clearly
imply $d\geq 2$, i.e. 
$t\geq 0$, hence $H^d_{I_{j_0}+\dots+I_{j_{t}}}(M)$ makes sense.

\medskip

\emph{Proof.} $c(I_jB)=c(I_j\hat A)\leq c$ by Lemma \ref{c(I)}. If
$p<t=[\frac{d-2}{c}]$, then
$d>c(p+1)$ and Theorem
\ref{main} with $\mathfrak p=p$ shows that
$H^q_{I_{j_0}+\dots+I_{j_p}}(M)=0$ for all $q>d-t+p$. Now Corollary
\ref{plan} shows that $H^{d-t}_I(M)=E^{-t,d}_2$ and Remark \ref{(-t,d)}
completes the proof. \qed

\medskip

This theorem implies some corollaries for arbitrary commmutative
Noetherian local rings containing a field. 
The following corollary deals with the vanishing of
$H^{v+1}_I(M)$ and generalizes Theorem
\ref{HL} to the case that
$\sqrt{IB}$ is not necessarily prime, but is the intersection of a small
number (less than $\frac{d}{c}$, to be precise) of prime ideals. Corollary
\ref{maincoro} is a specialization of Corollary \ref{Coro-gener} to the
case that the ring $A$ is complete, regular and strictly Henselian.

\begin{corollary}\label{Coro-gener}
Let $A$ be a local ring containing a field, let $c$ be a positive
integer and
let $I$ be an ideal of $A$ with $c(I\hat A)\leq c$. Let
$B$ be the completion of the strict Henselization of the
completion of $A$ and let $d>c$ be an integer. Assume that the ideal $IB$
has
$n<\frac{d}{c}$ minimal primes. 
If
$M$ is an
$A$-module with ${\rm dimSupp}M\leq d$, then
$H^q_I(M)=0$ for
$q>d-1-[\frac{d-2}{c}]$. 
In particular, if ${\rm dim}A\leq d$, then ${\rm cd}(A,I)\leq
d-1-[\frac{d-2}{c}]$.
\end{corollary}
\emph{Proof.} Let $v=d-1-[\frac{d-2}{c}]$. Theorem \ref{Falt} implies
that $H^q_I(M)=0$ for all $q>v+1$. If $c|(d-1)$, then
$v+1=d-[\frac{d-1}{c}]+1$, hence Theorem \ref{Falt} implies that
$H^{v+1}_I(M)=0$ in this case. Hence we only need to show that
$H^{v+1}_I(M)=0$ in the case that $c\not|(d-1)$. 

Since $B$ is faithfully
flat over
$A$ and 
$B\otimes_AH^{v+1}_I(M)\cong
H^{v+1}_{IB}(B\otimes_AM)$, 
it is enough to show that $H^{v+1}_{IB}(B\otimes_AM)=0$. The module
$H^{v+1}_{IB}(B\otimes_AM)$ is isomorphic to the cokernel of the
map $\Phi_{B\otimes_AM, IB}$ of Theorem \ref{general}. The facts that
$c\not|(d-1)$ and $n<\frac{d}{c}$ imply the bound $n\leq
[\frac{d-2}{c}]=t$. Hence
$(t+1)$-tuples $\{j_0,\dots,j_t\}$ such that $1\leq
j_o<\dots<j_t\leq n$ do not exist. Hence the target module of the
map $\Phi_{B\otimes_AM, IB}$ vanishes. Therefore the cokernel also
vanishes. \qed

\medskip

The following two corollaries, \ref{Artinian} and \ref{Delta}, deal with
the structure of $H^{v+1}_I(M)$ when this module doesn't
necessarily vanish.

In general local cohomology modules of
finitely generated modules are not necessarily Artinian, even if they are
supported in dimension zero \cite[Section
3]{Hart}. But Corollary \ref{Artinian} says, in particular, that
$H^{v+1}_I(M)$ is Artinian for a finitely generated $M$.

\begin{corollary}\label{Artinian} 
Let $A$ be a local ring containing a field, let $c>0$ and $d>1$ be
integers, let $v=d-1-[\frac{d-2}{c}]$, let $I$ be an ideal of $A$ with
$c(I\hat A)\leq c$ and let
$M$ be an $A$-module with ${\rm dimSupp}M\leq d$. Then 

(a) ${\rm dimSupp}H_I^{v+1}(M)=0$. If $M$ is finitely generated, then
$H_I^{v+1}(M)$ is Artinian.

(b) Let $B=\widehat{\hat A^{\rm sh}}$ be the completion of the strict
Henselization of the completion of $A$, let $k$ be a coefficient field of
$\hat A$ and let $K$ be the coefficient field of $B$ containing
$k$. Then
$H^{v+1}_{IB}(B\otimes_AM)=K\otimes_kH^{v+1}_I(M)$.
\end{corollary}

We note that every coefficient field $k$ of $\hat A$ can be
uniquely extended to a coefficient field $K$ of $B$ and $K$ is the
separable closure of
$k$ (in fact, $B\cong K\hat\otimes_k\hat A$ where $\hat \otimes_k$ is the
complete tensor product over $k$).

\medskip

\emph{Proof of Corollary \ref{Artinian}.} ${\rm dimSupp}(B\otimes_AN)={\rm
dimSupp}(N)$ for every $A$-module $N$ since
$B$ is faithfully flat over
$A$. Thus ${\rm dimSupp}(B\otimes_AM)\leq d$. This implies that 
each of the modules $H^d_J(B\otimes_AM)$ appearing in the
map $\Phi_{B\otimes_AM, IB}$ of Theorem \ref{general} is supported in
dimension zero. Indeed, if
$P$ is a non-maximal prime ideal of $B$, then ${\rm dim
Supp}(B\otimes_AM)_P<d$, hence
$H^d_J(B\otimes_AM)_P=H^d_{J_P}((B\otimes_AM)_P)$ vanishes at $P$. 

$c(IB)=c(I\hat A)\leq c$ by Lemma \ref{c(I)},
hence Theorem \ref{general} implies that $H^{v+1}_{IB}(B\otimes_AM)$ is
the cokernel of the map
$\Phi_{B\otimes_AM, IB}$. 
Thus 
$H^{v+1}_{IB}(B\otimes_AM)$ is the cokernel of
two modules supported in dimension zero and is therefore itself
supported in dimension zero. But $B\otimes_AH^{v+1}_I(M)\cong
H^{v+1}_{IB}(B\otimes_AM)$ since $B$ is faithfully flat over $A$. Hence
$H_I^{v+1}(M)$ also is supported in dimension zero. 

Let $N$ be a
finiely generated submodule of $H_I^{v+1}(M)$. Then $N$ is annihilated by
$\mathfrak m_A^s$ for some $s$, where $\mathfrak m_A$ is the maximal
ideal of $A$. The fact that $B/\mathfrak m_A^sB\cong K\otimes_kA/\mathfrak
m_A^s$ implies that $B\otimes_AN=K\otimes_kN$. Since $H_I^{v+1}(M)$ is
the direct limit of its finitely generated submodules and the tensor
product commutes with direct limits, $B\otimes_AH^{v+1}_I(M)\cong
K\otimes_kH^{v+1}_I(M)$.

Now assume $M$ is finitely generated. An $A$-module $N$ supported in
dimension zero is Artinian if and only if the $B$-module $K\otimes_kN$ is
Artinian. Hence it is enough to prove that
$H^{v+1}_{IB}(B\otimes_AM)=K\otimes_kH^{v+1}_I(M)$ is Artinian. Since
$H^{v+1}_{IB}(B\otimes_AM)$ is the cokernel of a map of two modules
each of which is a direct sum of modules of the form
$H^d_J(B\otimes_AM)$, it is enough to prove that each such module is
Artinian. 

$B\otimes_AM$ is a finitely generated
$B$-module. Hence there is a filtration $0=N_0\subset N_1\subset
N_2\subset\dots\subset N_s=B\otimes_AM$ such that $N_i/N_{i-1}\cong
B/P_i$ where the $P_i$s are some prime ideals of $B$. The resulting exact
sequences $H_J^d(N_{i-1})\to H_J^d(N_i)\to H_J^d(N_i/N_{i-1})$
imply by induction on $i$ that it is enough to prove that
$H^d_J(N_i/N_{i-1})=H^d_J(B/P_i)$ is Artinian for every $i$. By the
Hartshorne-Lichtenbaum local vanishing theorem \cite[3.1]{Ha},
$H^d_J(B/P_i)\ne 0$ only if
$\sqrt{J+P_i}=\mathfrak m_B$ (where $\mathfrak m_B$ is the maximal ideal
of
$B$) in which case
$H^d_J(B/P_i)=H^d_{\mathfrak m}(B/P_i)$. But $H^j_{\mathfrak m}(\mathcal
N)$ is Artinian for every $j$ and every finitely generated $B$-module
$\mathcal N$.
\qed

\medskip

If $M$ is finitely generated, Theorem \ref{general} dispels a
great deal of mystery about the module 
$H^{v+1}_{IB}(B\otimes_AM)$ by presenting it as the cokernel of an
explicit map between two Artinian $B$-modules, and then
the isomorphism
$H^{v+1}_{IB}(B\otimes_AM)=K\otimes_kH^{v+1}_I(M)$ of Corollary
\ref{Artinian} enables us to understand the $A$-module
$H^{v+1}_I(M)$ in terms of the $B$-module $H^{v+1}_{IB}(B\otimes_AM)$. In
particular, the length of the annihilator of
$\mathfrak m_A^s$ in
$H^{v+1}_I(M)$ equals the length of the annihilator of $\mathfrak m_B^s$
in $H^{v+1}_{IB}(B\otimes_AM)=K\otimes_kH^{v+1}_I(M)$ for every $s$.

\medskip

The following corollary shows that in
the case that
$M$ is supported on
$V(P)$ where $P$ is an ideal of $A$ such that 
$\sqrt{PB}$ is a prime ideal of $B$, the module $H^{v+1}_I(M)$ can be
completely described in terms of the module $H^d_{\mathfrak m}(M)$ and the
singular homology of a suitable simplicial complex. Theorem
\ref{maintheorem} is a specialization of Corollary \ref{Delta} to the
case that the ring $A$ is complete, regular and strictly Henselian.

\begin{corollary}\label{Delta}
Let $A$ be a local ring containing a field. Let $\mathfrak m$ and $k$ be
the maximal ideal and the residue field of $A$. Let
$c>0$ and
$d>1$ be integers, let
$I$ be an ideal of $A$ with $c(I\hat A)\leq c$ and let $B$ be the
completion of the strict Henselization of the completion of $A$.
Let $I_1, I_2,\dots, I_n$ be the minimal primes of $IB$, 
let $t=[(d-2)/c]$ and let
$v=d-1-[(d-2)/c]$. Let $P$ be an ideal of $A$ such that ${\rm
dim}(A/P)=d$ and the ideal $PB$ of $B$ has just one minimal
prime. Let $M$ be an
$A$-module supported on $V(P)$. Let
$\Delta$ be the simplicial complex on $n$ vertices $\{1,2,\dots,n\}$
defined as follows: a simplex
$\{j_0,\dots,j_s\}$ belongs to $\Delta$ iff
$I_{j_0}+\dots+I_{j_s}+PB$ is not $\mathfrak m_B$-primary, where
$\mathfrak m_B$ is the maximal ideal of $B$. Let
$w$ be the dimension over $k$ of $\tilde H_{t -1}(\Delta;k)$, the
$(t -1)$-th reduced singular homology group of $\Delta$ with
coefficients in
$k$. Then
$H^{v+1}_I(M)$ is isomorphic to the direct sum of $w$ copies of
$H^d_{\mathfrak m}(M)$. In particular, $H^{v+1}_I(M)=0$ if and only if
either
$H^d_{\mathfrak m}(M)=0$ or 
\hbox{$\tilde H_{t -1}(\Delta;k)=0$.}
\end{corollary}

\emph{Proof.} If $N$ is an $A$-module supported on $\{\mathfrak
m\}$, then $\hat A\otimes_AN=N$. This implies that 
$H^{v+1}_{I\hat A}(\hat A\otimes_AM)=\hat
A\otimes_AH^{v+1}_I(M)=H^{v+1}_I(M)$ because 
$H^{v+1}_I(M)$ is supported on $\{\mathfrak m\}$ by Corollary
\ref{Artinian} and $\hat A$ is flat over $A$. Hence we can replace $A$ by
$\hat A$ and $I$ by $I\hat A$, i.e. we can assume that $A$ is complete.
For the rest of the proof we assume that $A$ is complete. 

Let $k\subset A$ be a coefficient field of
$A$. Then $B=K\hat\otimes_kA$ where $K$ is the separable closure of
$k$. According to
\cite[4.2]{HL} there exists a finite Galois field extension
$\tilde k\supset k$ (of course $K\supset \tilde k$) such that upon setting
$\tilde A=\tilde k\otimes_kA$, the minimal primes $\tilde I_1,\dots,
\tilde I_n$ of
$\tilde I=I\tilde A$ are in a one-to one correspondence with the minimal
primes of $IB$, namely, $\sqrt{\tilde I_jB}=I_j$ for all $j$ (the ring
$\tilde A$ sits between $A$ and $B$). Clearly,
$\tilde A$ is a finite free (hence faithfully flat) $A$-module and a
complete local ring with residue field $\tilde k$ which it contains.
Since $B$ is the strict Henselization of $\tilde A$, it follows
from Lemma \ref{c(I)} that $c(\tilde I)=c(\tilde IB=IB)=c(I)$, hence
$c(\tilde I)\leq c$.

If $\mathcal N$ is an $\tilde A$-module and $\mathcal I$ is an ideal of
$A$, there is an isomorphism of $A$-modules $_AH^i_{\mathcal I\tilde
A}(\mathcal N)\cong H^i_{\mathcal I}(_A\mathcal N)$, where $_A(-)$ means
that the corresponding $\tilde A$-module is viewed as an $A$-module via
restriction of scalars. Accordingly, in the sequel we often denote
$H^i_{\mathcal I\tilde A}(\mathcal N)$ by 
$H^i_{\mathcal I}(\mathcal N)$ and view it as an $\tilde A$-module.

For
any
$A$-module
$N$ we set
$\tilde N=\tilde A\otimes_AN=\tilde k\otimes_kN$. We view $\tilde N$
as an $\tilde A$-module in a natural way. Since
$\tilde A$ is flat over $A$, there is an isomorphism of $\tilde A$-modules
$H^i_{\mathcal I}(\tilde N)\cong
\widetilde{H^i_{\mathcal I}(N)}$ for any ideal $\mathcal I$ of $A$.

By Theorem \ref{main}, if $p<t$, then 
$H^q_{\tilde I_{i_0}+\dots+\tilde I_{i_p}}(\tilde M)=0$ for all  
$q\geq d-t+p$ and all $\{i_0,\dots,i_p\}$. Hence by Corollary
\ref{plan}, $H^{d-t=v+1}_{I\tilde A}(\tilde M)\cong
E_2^{-t,d}$. Clearly,
$\tilde M$ is supported on
$V(\tilde P)$ where $\tilde P=\sqrt{P\tilde A}$ is a prime ideal of
$\tilde A$ such that dim$(\tilde A/\tilde P)=d$ and the ideal $\tilde
PB=PB$ has just one minimal prime. Hence Corollary \ref{MVDelta} implies
 $E_2^{-t,d}\cong\tilde
H_{t-1}(\Delta;H^d_{\mathfrak m}(\tilde M))\cong  H^d_{\mathfrak
m}(\tilde M)\otimes_{\tilde k}\tilde H_{t-1}(\Delta;\tilde k)$ 
where the last isomorphism holds by the universal
coefficients theorem since we are over a field. But dim$_{\tilde k}\tilde
H_{t-1}(\Delta;\tilde k)={\rm dim}_k\tilde H_{t-1}(\Delta;k)=w$ because
the field $\tilde k$ has the same characteristic as the field $k$.
Hence the module 
$H^{v+1}_I(\tilde M)\cong\widetilde{H^{v+1}_I(M)}$ is isomorphic to the
direct sum of
$w$ copies of the module $H^d_{\mathfrak m}(\tilde M)\cong
\widetilde{H^d_{\mathfrak m}(M)}$, i.e. $\widetilde{H^{v+1}_I(M)}\cong
\widetilde{H^d_{\mathfrak m}(M)^w}$ in the category of
$\tilde A$-modules. This completes the proof in the
special case that $A=\tilde A$, i.e. every minimal prime of $IB$ is the
radical of the extension to $B$ of a minimal prime of $I$. The rest of
the proof consists in deducing the general case from this special case.
This deduction turns out to be an unexpectedly long story.

Let $G$ be the Galois group of $\tilde k$ over $k$. Let
$A[G]=\{\sum_{\sigma\in G} a_{\sigma}\sigma|a_{\sigma}\in A\}$ be the
group ring of $G$ with coefficients in $A$; every $a\in A$ commutes with
every
$\sigma\in G$. By an $A[G]$-module we always mean a left $A[G]$-module.
If $N$ is an $A$-module, $\tilde N=\tilde k\otimes_kN$
acquires a {\it standard structure} of $A[G]$-module via $\sigma(c\otimes
x)=\sigma(c)\otimes x$ for all $\sigma\in G, c\in \tilde k$ and $x\in N$.
In fact $(\tilde-)$ is a functor from $A$-mod to $A[G]$-mod, namely, if
$\phi:N\to N'$ is a morphism of $A$-modules, then $\tilde\phi={\rm
id}\otimes\phi:\tilde k\otimes_kN\to \tilde k\otimes_kN'$ is a morphism
of $A[G]$-modules. 

Let
$J\subset A[G]$ be the ideal generated by the set
$\{\sigma -\tau|\sigma, \tau\in G\}$. \hbox{Then $N\cong \tilde N/J\tilde
N$ for every $A$-module $N$.}
To complete the proof in the general case it is enough to show the
existence of an
$A[G]$-module isomorphism $\widetilde{H^{v+1}_I(M)}\cong
\widetilde{H^d_{\mathfrak m}(M)^w}$ where the
$A[G]$-module structures on $\widetilde{H^{v+1}_I(M)}$ and 
$\widetilde{H^d_{\mathfrak m}(M)}$ are the standard ones (it has been
shown above that the two modules are isomorphic as $\tilde
A$-modules). Indeed, 
$\widetilde{H^{v+1}_I(M)}/J\widetilde{H^{v+1}_I(M)}\cong H^{v+1}_I(M)$
and $\widetilde{H^d_{\mathfrak m}(M)}/J\widetilde{H^d_{\mathfrak
m}(M)}\cong H^d_{\mathfrak m}(M)$, hence the existence of the
above-mentioned $A[G]$-module isomorphism would imply that
$H^{v+1}_I(M)\cong H^d_{\mathfrak m}(M)^w$, as needed. The rest of the
proof is devoted to showing the existence of such an $A[G]$-module
isomorphism.

If
$N$ is an $A[G]$-module, then $H^i_I(N)$ inherits a
structure of $A[G]$-module. Namely, the action $\sigma:N\to
N$, being an $A$-module homomorphism, induces the action
$\sigma:H^i_I(N)\to H^i_I(N)$ for every $\sigma\in G$.
Thus we get a functor $H^i_I(-):A[G]{\rm -mod}\to A[G]{\rm -mod}$. It is
easy to check that for every $A$-module $N$ there is an isomorphism of
$A[G]$-modules $H^i_I(\tilde N)\cong \widetilde{H^i_I(N)}$ where the
$A[G]$-module structures on $\tilde N$ and $\widetilde{H^i_I(N)}$ are the
standard ones and the $A[G]$-module structure on $H^i_I(\tilde N)$ is
functorially induced from $\tilde N$.

If $N$ is an $A$-module and $E_A(N)$ is the
injective hull of $N$ in the category of $A$-modules, then
$\widetilde{E_A(N)}\cong E_{\tilde A}(\tilde N)$ where $E_{\tilde
A}(\tilde N)$ is the injective hull of $\tilde N$ in the category of
$\tilde A$-modules, i.e.
$\widetilde{E_A(N)}$ is an injective $\tilde A$-module. Since every
injective $A$-module is a direct sum of modules of the form $E_A(N)$, we
conclude that if $E$ is an injective $A$-module, then $\tilde E$ is an
injective $\tilde A$-module. Hence if the complex $E^{\bullet}$ is an
injective resolution of $M$ in the
category of
$A$-modules, then $\widetilde {E^{\bullet}}$ is both an injective
resolution of
$\tilde M$ in the category of $\tilde A$-modules and a complex in the
category of $A[G]$-modules.

The group $G$ acts on $\tilde A$ in a standard way (i.e. $\sigma(c\otimes
a)=\sigma(c)\otimes a$ for all $\sigma\in G, c\in \tilde k$ and $a\in
A$) and the elements of
$G$ permute the ideals
$I_0\dots,I_n$. For any
$\sigma\in G$ we write $\sigma(i)=j$ if $\sigma(I_i)=I_j$. We recall
from the proof of Theorem \ref{MV} that $\Gamma^{-p}(\tilde M)=\oplus_{
i_0<i_1<\dots<i_p}\Gamma_{I_{i_0}+\dots+I_{i_p}}(\tilde M)$. This is an
$A[G]$-module, for if $x\in \Gamma_{I_{i_0}+\dots+I_{i_p}}(\tilde M)$,
then $\sigma(x)\in \Gamma_{I_{\sigma(i_0)}+\dots+I_{\sigma(i_p)}}(\tilde
M)$ (indeed, if $b\in \tilde
A, x\in \tilde M$ and $\sigma\in G$, then $\sigma(bx)=\sigma(b)\sigma(x)$,
hence an ideal $\mathcal I$ of $\tilde A$ annihilates $x$ if and only if
the ideal $\sigma(\mathcal I)$ annihilates $\sigma(x)$). 
It is not hard to see that the differentials in
the complex
$\Gamma^{\bullet}(\tilde M)$ are $A[G]$-module homomorphisms. Hence the
double complex
$\Gamma^{\bullet}(\Gamma_I(\widetilde{E^{\bullet}}))$ that induces the
Mayer-Vietoris spectral sequence 
$$E_1^{-p,q}=\oplus_{i_0<\dots<i_p}H^q_{I_{i_0}+\dots+I_{i_p}}(\tilde M)
\Longrightarrow H^{q-p}_I(\tilde M)\cong\widetilde{H^{q-p}_I(M)}$$
(see the proof of Theorem \ref{MV}) is
a double complex in the category of $A[G]$-modules. Hence this
 is a spectral
sequence in the category of $A[G]$-modules.

The map $\Phi_{\tilde M,\tilde I}$ of Theorem \ref{general} is a map on
the
$E_1$ page of the spectral sequence, hence it is an $A[G]$-module map. The
same argument as in the proof of Theorem
\ref{general} shows that the abutment of the spectral sequence in degree
$v+1$ is isomorphic (in the category of $A[G]$-modules) to the cokernel of
$\Phi_{\tilde M,\tilde I}$. But the abutment is $\widetilde{H^{v+1}_I(
M)}$ with the standard $A[G]$-module structure. Hence it is enough to
show that the cokernel of $\Phi_{\tilde M,\tilde I}$ is isomorphic as an
$A[G]$-module to $\widetilde{H^{v+1}_I(M)^w}$.

For an integer $s$ let $\Lambda_s$ be the set of the subsets
$\{j_0,\dots, j_s\}\subset \{1,\dots,n\}$ of cardinality $s+1$ 
such that $\sqrt{\tilde I_{j_0}+\dots+\tilde I_{j_s}+\tilde P}=\tilde {\mathfrak m}$. The Hartshorne-Lichtenbaum local vanishing theorem
\cite[3.1]{Ha} implies that
$H^d_{I_{j_0}+\dots+I_{j_s}}(\tilde M)$ vanishes if 
$\{j_0,\dots, j_s\}\not\in \Lambda_s$, otherwise
$H^d_{I_{j_0}+\dots+I_{j_s}}(\tilde M)\cong H^d_{\mathfrak m}(\tilde M)$.
For
$\lambda=\{j_0,\dots,j_s\}\in
\Lambda_s$ we denote $H^d_{I_{j_0}+\dots+I_{j_s}}(\tilde M)$ 
by $H^d_{\mathfrak
m}(\tilde M)_{\lambda}.$ The module $H^d_{\mathfrak
m}(\tilde M)_{\lambda}$ is just a copy of 
$H^d_{\mathfrak m}(\tilde M)$ indexed by $\lambda$. The map
$\Phi_{\tilde M, I}$ of Theorem
\ref{general} takes the form 
$$\Phi:\oplus_{\lambda\in\Lambda_{t+1}}H^d_{\mathfrak
m}(\tilde M)_{\lambda}
\to \oplus_{\lambda\in\Lambda_{t}}H^d_{\mathfrak
m}(\tilde M)_{\lambda}.$$
It follows from the description of the map $\Phi_{M,I}$ in the statement
of Theorem \ref{general} that if $\lambda'\in \Lambda_{t+1}$ and
$\lambda\in \Lambda_t$, then the map
$H^d_{\mathfrak m}(M)_{\lambda'}\to H^d_{\mathfrak m}(M)_{\lambda}$
induced by $\Phi$ is
$0$ if $\lambda\not\subset\lambda'$ and is either $\rm id$ or $-{\rm id}$
if $\lambda\subset\lambda'$ where $\rm id$ denotes the identity map on
$H^d_{\mathfrak m}(M)$ (we view $H^d_{\mathfrak m}(M)_{\lambda'}$ and 
$H^d_{\mathfrak m}(M)_{\lambda}$ as two different copies of the same
module
$H^d_{\mathfrak m}(M)$).

Let $y\in \tilde k$ be an element such that the set $\{\sigma(y)|\sigma\in
G\}$ is a $k$-basis of $\tilde k$. There is an isomorphism of $A$-modules
$\tilde M\cong \oplus_{\sigma\in G}(\sigma(y)\otimes_kM)$ and the fact
that
$G$ acts on $\tilde M$ in a standard way means that 
$\sigma'(\sigma(y)\otimes_kx)=\sigma'\sigma(y)\otimes_kx$ for all $x\in
M$ and all $\sigma',\sigma\in G$. Hence $H^d_{\mathfrak
m}(\tilde M)_{\lambda}\cong\oplus_{\sigma\in
G}(\sigma(y)\otimes_kH^d_{\mathfrak m}(M)_{\lambda})$ (the only function
of the subscript
$\lambda$ in $\sigma(y)\otimes_kH^d_{\mathfrak m}(M)_{\lambda}$ is 
to indicate that this particular copy of
$\sigma(y)\otimes_kH^d_{\mathfrak m}(M)$ came from 
$H^d_{\mathfrak m}(\tilde M)_{\lambda}$). Hence the map $\Phi$ takes the
following form:
$$\Phi:\oplus_{\sigma\in G, \lambda\in
\Lambda_{t+1}}(\sigma(y)\otimes_kH^d_{\mathfrak m}(M)_{\lambda})\to 
\oplus_{\sigma\in G, \lambda\in
\Lambda_{t}}(\sigma(y)\otimes_kH^d_{\mathfrak m}(M)_{\lambda}).$$
The map $\sigma'(y)\otimes_kH^d_{\mathfrak m}(M)_{\lambda'}\to
\sigma(y)\otimes_kH^d_{\mathfrak m}(M)_{\lambda}$ induced by $\Phi$ is
$0$ if $\sigma'\ne \sigma$ or $\lambda\not\subset\lambda'$ and is
either $\rm id$ or $-{\rm id}$ if $\sigma'=\sigma$ and
$\lambda\subset\lambda'$ where $\rm id$ denotes the identity map on
$\sigma(y)\otimes_kH^d_{\mathfrak m}(M)$. 

The above description of the $G$-action on $\Gamma^{-p}(\tilde M)$
implies that $G$ acts on $\oplus_{\sigma\in G, \lambda\in
\Lambda_{s}}(\sigma(y)\otimes_kH^d_{\mathfrak m}(M)_{\lambda})$, where
$s=t, t+1$, as follows. An element
$\sigma'\in G$ sends the $A$-module
$\sigma(y)\otimes_kH^d_{\mathfrak m}(M)_{\lambda}$ isomorphically onto
the $A$-module $\sigma'\sigma(y)\otimes_kH^d_{\mathfrak
m}(M)_{\sigma'(\lambda)}$ via the map $\sigma(y)\otimes x\mapsto
\sigma'\sigma(y)\otimes_kx$ for every
$x\in H^d_{\mathfrak m}(M)$ (if $\lambda=\{j_0,\dots,j_s\}$, then
$\sigma'(\lambda)=\{\sigma'(j_0),\dots, \sigma'(j_s)\}$). 

Let $O_t$ (resp. $O_{t+1}$) be the set of the orbits of the action of
$G$ on the set of the $A$-modules
$\{\sigma(y)\otimes_kH^d_{\mathfrak m}(M)_{\lambda} | \sigma\in
G, \lambda\in \Lambda_t$ ({\rm resp.}\ 
$\lambda\in \Lambda_{t+1}$)\}. Since $\sigma'\sigma(y)=\sigma(y)$ if and
only if $\sigma'$ is the identity element of $G$, the stabilizer of each
module $\sigma(y)\otimes_kH^d_{\mathfrak m}(M)_{\lambda}$ is trivial.
Hence each orbit $\mathfrak o$ consists of
$|G|$ elements which are $\{\sigma(y)\otimes_kH^d_{\mathfrak
m}(M)_{\sigma(\lambda)}|\sigma\in G\}$, where $\lambda$ is fixed.
We denote by $\widetilde{H^d_{\mathfrak m}(
M)}_{\mathfrak o}$ the direct sum of the $|G|$ elements of
$\mathfrak o$, i.e. $\widetilde{H^d_{\mathfrak m}(
M)}_{\mathfrak o}=\oplus_{\sigma\in G}(\sigma(y)\otimes_kH^d_{\mathfrak
m}(M)_{\sigma(\lambda)})$. Ignoring the subscripts $\sigma(\lambda)$ and
$\mathfrak o$ (which are there just to keep track of
different copies of the same modules)
 this direct sum is isomorphic to the
$A[G]$-module
$\widetilde{H^d_{\mathfrak m}(M)}$ with the standard $A[G]$-module
structure. 

The above descriptions of the map $\Phi$ and the $G$-action on
its source and target, and the fact that $\Phi$ commutes with this
$G$-action, imply that the map
$\sigma'(y)\otimes_kH^d_{\mathfrak m}(M)_{\sigma'(\lambda')}\to
\sigma(y)\otimes_kH^d_{\mathfrak m}(M)_{\sigma(\lambda)}$ induced by
$\Phi$ is
$0$ if
$\sigma'\ne\sigma$, and if $\sigma'=\sigma$, then this map 
$\sigma(y)\otimes_kH^d_{\mathfrak m}(M)_{\sigma(\lambda')}\to
\sigma(y)\otimes_kH^d_{\mathfrak m}(M)_{\sigma(\lambda)}$ 
depends only on $\lambda'$ and $\lambda$ but
not on $\sigma$, and is either 
$0$, or
$\rm id$, or
$-{\rm id}$ where $\rm id$ denotes the identity map on
$\sigma(y)\otimes_kH^d_{\mathfrak m}(M)$. Hence every map 
$\widetilde{H^d_{\mathfrak m}(
M)}_{\mathfrak o'}\to \widetilde{H^d_{\mathfrak m}(
M)}_{\mathfrak o}$ induced by $\Phi$ is either $0$, or $\rm id$, or
$-{\rm id}$ where $\rm id$ denotes the identity map on 
$\widetilde{H^d_{\mathfrak m}(
M)}$.

Hence the map $\Phi$ takes the form 
$$\Phi:\oplus_{\mathfrak o\in O_{t+1}}\widetilde{H^d_{\mathfrak
m}(M)}_{\mathfrak o}\to\oplus_{\mathfrak o\in O_t}
\widetilde{H^d_{\mathfrak m}(M)}_{\mathfrak o}$$
where each $\widetilde{H^d_{\mathfrak m}(M)}_{\mathfrak o}$ is isomorphic
to the same $A[G]$-module $\widetilde{H^d_{\mathfrak m}(M)}$
and the matrix defining $\Phi$ is an
$|O_{t+1}|\times|O_t|$-matrix
$M$ each of whose entries is either 0, or 1, or -1. Since the entries of
$M$ are all in $\tilde k$ and $\tilde k$ is a field, the cokernel of
$\Phi$ is isomorphic as an
$A[G]$-module to the direct sum of $w'$ copies of
$\widetilde{H^d_{\mathfrak m}(M)}$ where $w'$ is the dimension over
$\tilde k$ of the cokernel of the map $\tilde k^{|O_{t+1}|}\to
\tilde k^{|O_t|}$ defined by the same matrix $M$. 

It has been shown earlier that the cokernel of $\Phi$ is isomorphic as an
$A$-module to a direct sum of $w$ copies of $\widetilde{H^d_{\mathfrak
m}(M)}$. Hence there is an isomorphism of $A$-modules
$\widetilde{H^d_{\mathfrak m}(M)^w}\cong \widetilde{H^d_{\mathfrak
m}(M)^{w'}}$. Let $\ell$ be the dimension over $k$ of the socle of
$\widetilde{H^d_{\mathfrak m}(M)}$. Since $\widetilde{H^d_{\mathfrak
m}(M)}$ is supported at $\mathfrak m$, $\ell=0$ if and only if $\widetilde{H^d_{\mathfrak
m}(M)}$. Hence $\ell=0$ implies $H^d_{\mathfrak m}(M)=0$ because
$\tilde A$ is faithfully flat over $A$. Hence if $\ell=0$, there is
nothing to prove. Now we assume
$\ell\ne 0$. The dimensions of the socle of the
$A$-module
$\widetilde{H^d_{\mathfrak m}(M)^w}\cong \widetilde{H^d_{\mathfrak
m}(M)^{w'}}$ is, on the one hand, $\ell w$ and on the other, $\ell w'$.
The equality $\ell w=\ell w'$ implies $w=w'$ since $\ell\ne0$. Hence the
cokernel of
$\Phi$ is isomorphic as an $A[G]$-module to 
$\widetilde{H^d_{\mathfrak m}(M)^w}$.      \qed 

\medskip

The result of the preceding corollary is based on the fact that if $M$
is supported on $V(P)$ where $P$ is an ideal of $A$ such that $\sqrt{PB}$
has just one minimal prime, then the set of the isomorphism classes of
the modules $H^d_{I_{j_0}+\dots+I_{j_{s}}}(\tilde M)$ appearing in the
map $\Phi_{\tilde M,\tilde I}$ has at most two
elements, 0 and $H^d_{\mathfrak m}(\tilde M)$. But if
$M$ is supported on
$V(P)$ where
$P\subset A$ is such that the ideal 
$PB$ has more than one minimal primes, then the set of the isomorphism
classes of the modules $H^d_{I_{j_0}+\dots+I_{j_{s}}}(\tilde M)$ may be
quite big. This precludes expressing the module $H^{v+1}_I(M)$ in terms
of a single module $H^d_{\mathfrak m}(M)$ and some singular homology. But
one can still give a necessary and sufficient
condition for cd$(A,I)\leq v$ in terms of singular homology alone, as
our final corollary shows.  
\begin{corollary}
Let $A$ be a $d$-dimensional local ring containing a field. 
Let $\mathfrak m$ and $k$ be
the maximal ideal and the residue field of $A$. Let
$c>0$ and $d>1$ be integers, let $t=[(d-2)/c]$ and let
$v=d-1-[(d-2)/c]$. Let
$I$ be an ideal of $A$ such that $c(I\hat A)\leq c$. Let $B$ be
the completion of the strict Henselization of the completion of
$A$. Let $I_1, I_2,\dots I_n$ be the minimal primes of $IB$ and let
$P_1,P_2,\dots$ be the minimal
$d$-dimensional primes of $B$ (a prime $P$ of $B$ is
$d$-dimensional if $B/P=d$). Let
$\Delta_i$ be the simplicial complex on $n$ vertices $\{1,2,\dots,n\}$
such that a simplex
$\{j_0,\dots,j_s\}$ belongs to $\Delta_i$ if and only if 
$I_{j_0}+\dots+I_{j_s}+P_i$ is not $\mathfrak
mB$-primary. Let $\tilde H_{t -1}(\Delta_i;k)$ be the
$(t -1)$-th reduced singular homology group of $\Delta_i$ with
coefficients in
$k$. Then ${\rm cd}(A,I)\leq v$ if and only if 
$\tilde H_{t -1}(\Delta_i;k)=0$ for every $i$.
\end{corollary}

\emph{Proof.} If $M$ is a $B$-module, then
$H^{v+1}_I(_AM)\cong _AH^{v+1}_{IB}(M)$ where the subscript $_A-$ means
that the corresponding $B$-module is viewed as an $A$-module via
restriction of scalars. If
$M$ is an
$A$-module, $H^{v+1}_{IB}(B\otimes_AM)\cong
B\otimes_AH^{v+1}_I(M)$ vanishes if and only if
$H^{v+1}_I(M)$ vanishes since $B$ is faithfully flat over $A$. This
implies that cd$(A,I)={\rm cd}(B,IB)$. 

If $\tilde H_{t -1}(\Delta_i;k)\ne 0$ for
some
$i$, then $H^{v+1}_I(B/P_i)\ne 0$ by Corollary \ref{Delta} since
$H^d_{\mathfrak mB}(B/P_i)\ne 0$. Conversely, assume $\tilde H_{t
-1}(\Delta_i;k)=0$ for every
$i$. We need to show that $H^{v+1}_I(M)=0$ for every $B$-module $M$.
Since every module is the direct limit of its finitely
generated submodules and $H^{v+1}_I(-)$ commutes with direct limits, we
can assume that $M$ is finitely generated in which case there is a finite
filtration
$0= M_0\subset M_1\subset\dots \subset M_{s-1}\subset M_s=M$ such that
$M_i/M_{i-1}\cong B/Q_i$ where $Q_i$ is some prime ideal of $B$.
Applying the functor $H^{v+1}_I(-)$ to the exact sequences $0\to
M_{i-1}\to M_i\to M_i/M_{i-1}\to 0$ and using induction on $i$, it is
enough to prove that
$H^{v+1}_I(B/Q)=0$ for every prime ideal $Q$ of $B$. 

If dim$B/Q<d$, then
$H^{v+1}_I(B/Q)=0$ by Theorem \ref{Falt}. If dim$B/Q=d$, then $Q=P_i$ for
some $i$ in which case
$H^{v+1}_I(B/Q)=0$ by Corollary \ref{Delta} because $\tilde H_{t
-1}(\Delta_i;k)=0$. \qed

\end{document}